%% file: main.tex
\documentclass[review]{elsarticle}

\usepackage{lineno,hyperref}
\modulolinenumbers[5]
\usepackage{comment}
\journal{Journal of \LaTeX\ Templates}

\usepackage[english]{babel}	
\usepackage[usenames]{color}
\usepackage{algorithm}
\usepackage[noend]{algpseudocode}

\usepackage[T1]{fontenc}
\usepackage{tgbonum}
\usepackage{utopia}
\usepackage{tgadventor}
\usepackage{lmodern}
\usepackage{tgcursor}
\usepackage{palatino}

\usepackage{color}
\usepackage{colortbl}
\usepackage{amsmath}
\usepackage{float} 
\usepackage{cmap}
\DeclareUnicodeCharacter{2212}{-}

\begin{document}

\begin{frontmatter}






\title{Agent-based modeling of COVID-19 outbreaks for New York state, UK and Novosibirsk region}

\author[1,2]{Olga Krivorotko\corref{cor1}}
\cortext[cor1]{Corresponding author}
\ead{krivorotko.olya@mail.ru}

\author[2]{Mariia Sosnovskaia}
\ead{m.sosnovskaya@g.nsu.ru}

\author[2]{Ivan Vashchenko}
\ead{i.vashchenko@g.nsu.ru}

\author[3]{Cliff Kerr}
\ead{ckerr@idmod.org}

\author[4]{Daniel Lesnic}
\ead{D.Lesnic@leeds.ac.uk}

\address[1]{Institute of Computational Mathematics and Mathematical Geophysics Siberian Branch of the Russian Academy of Sciences, 6 Prospect Akademika Lavrentieva Street, Novosibirsk, 630090, Russia;}
\address[2]{Novosibirsk State University, 2 Pirogova Street, Novosibirsk, 630090, Russia;}
\address[3]{Institute for Disease Modeling, Bill \& Melinda Gates Foundation, Seattle, USA.}
\address[4]{University of Leeds, LS2 9JT, UK.}

\input{1abstract}
\newpageafter{1abstract} 
\end{frontmatter}

\linenumbers 

\input{2introduction}
%
\input{3data_analysis}
%
\input{4agent_model}
%
\input{5methods_and_approaches}
%
\input{6modelling_and_forecasting}

\input{7conclusion}

\input{9competing}

\input{acknowledges}

\input{references}

\end{document}

%% file: 1abstract.tex
\begin{abstract}
This paper uses Covasim, an agent-based model (ABM) of COVID-19, to evaluate and scenarios of epidemic spread in New York State (USA), the UK, and the Novosibirsk region (Russia). Epidemiological parameters such as contagiousness (virus transmission rate), initial number of infected people, and probability of being tested depend on the region’s demographic and geographical features, the containment measures introduced; they are calibrated to data about COVID-19 spread in the region of interest. At the first stage of our study, epidemiological data (numbers of people tested, diagnoses, critical cases, hospitalizations, and deaths) for each of the mentioned regions were analyzed. The data were characterized in terms of seasonality, stationarity, and dependency spaces, and were extrapolated using machine learning techniques to specify unknown epidemiological parameters of the model. At the second stage, the Optuna optimizer based on the tree Parzen estimation method for objective function minimization was applied to determine the model’s unknown parameters. The model was validated with the historical data of 2020. The modeled results of COVID-19 spread in New York State, the UK and the Novosibirsk region have demonstrated that if the level of testing and containment measures is preserved, the number of positive cases in New York State and the Novosibirsk region will remain the same during March of 2021, while in the UK it will reduce. Due to the features of the data for the Novosibirsk region (two datasets are characterized as stationary series with probability of 1), the forecast precision is relatively high  for the number of hospitalizations, but is lower for new cases of COVID-19. 

\end{abstract}


%% file: 2introduction.tex
\section{Introduction}

In December 2019, the COVID-19 pandemic originated in the Wuhan province of China. Since that time, more than 130 million people in 192 countries have been infected with the disease, and more than 2.87 million people have died after getting infected. During the year 2020, mankind mobilized its resources to fight the pandemic. One of the useful tools in this struggle has been mathematical modeling that uses known historical data to study different scenarios of disease spread \cite{1, 2}. 

The models including those for studying coronavirus infections can be divided into two groups: compartmental and agent-based models. In compartmental models, a population is divided into groups sharing similar features and interacting with one another following the mass action law. Agent-based models (ABMs), on the other hand, give each agent (people, social institutes, the state, etc.) a set of features and determines the way the agents interact from random graphs following disease spread principles. In other words, an agent’s behavior is determined individually, and their joint behavior is described as the interaction of multiple agents (bottom-top approach). Unlike compartmental models, ABMs are capable of providing a detailed description of an epidemiological situation, especially in a case of inaccurate and insufficient data. ABMs account for the stochastic nature of epidemic spread, makes it possible to estimate the likelihood of different epidemic scenarios, and allows one to evaluate the risks of unfavorable events occurring due to policy changes. The resulted data enable one to make conclusions about the duration, severity and scale of an epidemic, evaluate the efficacy of the preventive and quarantine measures, and assess its economic consequences.

During 2020, many models were developed to predict COVID-19 spread. In \cite{masking}, both compartmental and agent-based models are presented to study mass face-mask wearing and predict its effect on COVID-19 spread. A graph-based ABM is suggested in \cite{wolframgraph}. The paper considers a small population of 1000 agents and graphs of different kinds such as fully-connected, Barabasi-Albert, Watts-Strogatz ones, etc. Another approach to agents interactions is demonstrated in \cite{ABM_grid}. In this paper, agents are initialized in a 2D space (so-called mesh) to consider a distance between different agents, so such a factor as social distancing can be explicitly accounted for. However, this approach leads to excessive computation complexity, and for that reason, these models consider a small number of agents (up to 500).

Another important direction of modeling has been a comparison of containment measures and projecting the future for different scenarios. The paper published by Chang et al. \cite{Nature} combines an epidemiological SEIR model and the hourly GPS data from the mobile phones of 98 million people in 10 US cities. The model predicted that closing of the most crowded public places such as restaurants and religious establishments would be a sufficient measure to contain the pandemic unlike unilateral measures to limit people’s mobility. In \cite{ABM_econ_brazil}, the authors offer an ABM using SEIR agents to model COVID-19 spread dynamics, where the agents imitate people, businesses and the government. Using the model, they have analyzed seven social-distancing scenarios having different epidemiological and economic effects. The paper has demonstrated that the so-called vertical isolation has no positive effect. In \cite{ABM_econ_german}  macroeconomic epidemiological ABM been presented to study the economic effect COVID-19 would have in different scenarios of pandemic containment such as closing of educational and entertainment facilities. The model was calibrated using the statistical data on country’s and business demography, households, employment, profits and wages in Germany.  

It is noteworthy that data collecting and processing is a very important step in building an effective COVID-19 spread model. However, in the studies mentioned above data pre-processing for the modeled regions was not performed. Most of them concentrated on building the models and algorithms, whose parameters were considered known either from literature or from experts’ estimations, so the issues of identifiability for unknown parameters have remained unresolved, as has the need of devising a regularization algorithm for solving the problem of epidemiological forecasting.  

In this study, our focus is on the analysis of data, parameter identification, and regularization algorithms. However, it is a known fact that epidemics develop differently in different locations. To address this issue, in our study, the epidemiological situations in New York State (USA), the UK, and the Novosibirsk region (Russia) were compared and analyzed.  

\subsection{New York State}
Nearly 2 million people were confirmed infected and more than 50000 people died by April 6, 2021 in NY State \cite{TimesNY}. In NY State, the pandemic spread rapidly, reaching its peak in March-April of 2020 (see Section \ref{sect: NY}). The healthcare system was overloaded in the very first months of the outbreak. Thanks to the containment measures introduced, the number of infected people reduced to a characteristic plateau that was followed by a second infection wave several months later.  

\subsection{United Kingdom}

Another pandemic scenario was observed in the UK. After the first infection wave had been successfully suppressed by June, 2020, the second wave hit the country hard and, due to the B.1.1.7 SARS-CoV-2 variant, the number of infected people increased 10 times compared to the first wave (see Section  \ref{sect: UK}). At the time of writing, this has reduced to 7000 cases a day, and the pandemic is on wane. 

\subsection{Novosibirsk region}
The region’s peculiar feature was the absence of the characteristic peaks that rather looked like plateaus (see Section 2.3), so the lockdown introduced in April 2020 had outlived its usefulness by the summer of 2020.  
 
This study is organized as follows. First, open-source data of COVID-19 spread in New York State, the UK and the Novosibirsk region were processed and analyzed using the  statistics and machine learning methods to find interdependencies, study seasonality and predict possible future dynamics (see Section 2). Second, we confirmed that the selected ABM met the identifiability condition such as being sensitive to data errors and capable of unambiguous determination of the unknown parameters of COVID-19 spread from additional measurements (Section \ref{sect: inv problem}). The obtained space of identifiable parameters was specified using the multilevel global-optimization method (Section \ref{sect: methods}). Finally, scenarios of how the COVID-19 pandemic that could develop in New York State, the UK and the Novosibirsk region were assessed concerning available data and certain containment measures (Section \ref{sect: results}).

%% file: 3data_analysis.tex
\section{Data analysis}
Data analysis and data processing are important parts of forecasting modeling. Before data processing begins, one has to understand the character of available data and determine their features. Collecting such daily indicators as the number of tests, diagnosed cases, ventilated COVID patients, etc. (see Table~\ref{tab:data_statistics}) helps to overview the general picture for a considered region, while anomalous time intervals may call for more scrupulous analysis.
\begin{table}[H]
\begin{center}
\caption{Indicators used in data analyses and their description.}
\label{tab:data_statistics}
\begin{tabular}{|p{2.5cm}|p{9cm}| } 
\hline
Indicator & Description\\ 
\hline
New Tests & Number of performed tests\\
New Diagnoses & Number of diagnosed cases\\
New Deaths & Number of deaths related to a positive COVID-19 diagnosis \\
Num Critical & Number of ICU patients with suspected and positive COVID-19 diagnosis\\

\hline
\end{tabular}
\end{center}
\end{table}

\subsection{New York State (USA)} \label{sect: NY}
To forecast the way the pandemic would develop, the data from the COVID Tracking Project’s website  \cite{covid_tracking} were used. The site contains information for each state and for the country as a whole. The feature of the data in question is the method they calculate  {\fontfamily{cmss}\selectfont New Diagnoses}: positive cases (confirmed plus probable) summing the total number of confirmed cases and the probable cases of COVID-19 reported by the state or territory, ideally per the ”August 5, 2020 CSTE case definition”. Some states are following the older ”April 5th, 2020 CSTE case definition” or using their own custom definitions. The latter method of data collection is more suitable for ABMs since it takes into account the percentage of infected people that may have been simply neglected. 

\begin{figure}[H]
    \centering
    \includegraphics[width=\textwidth]{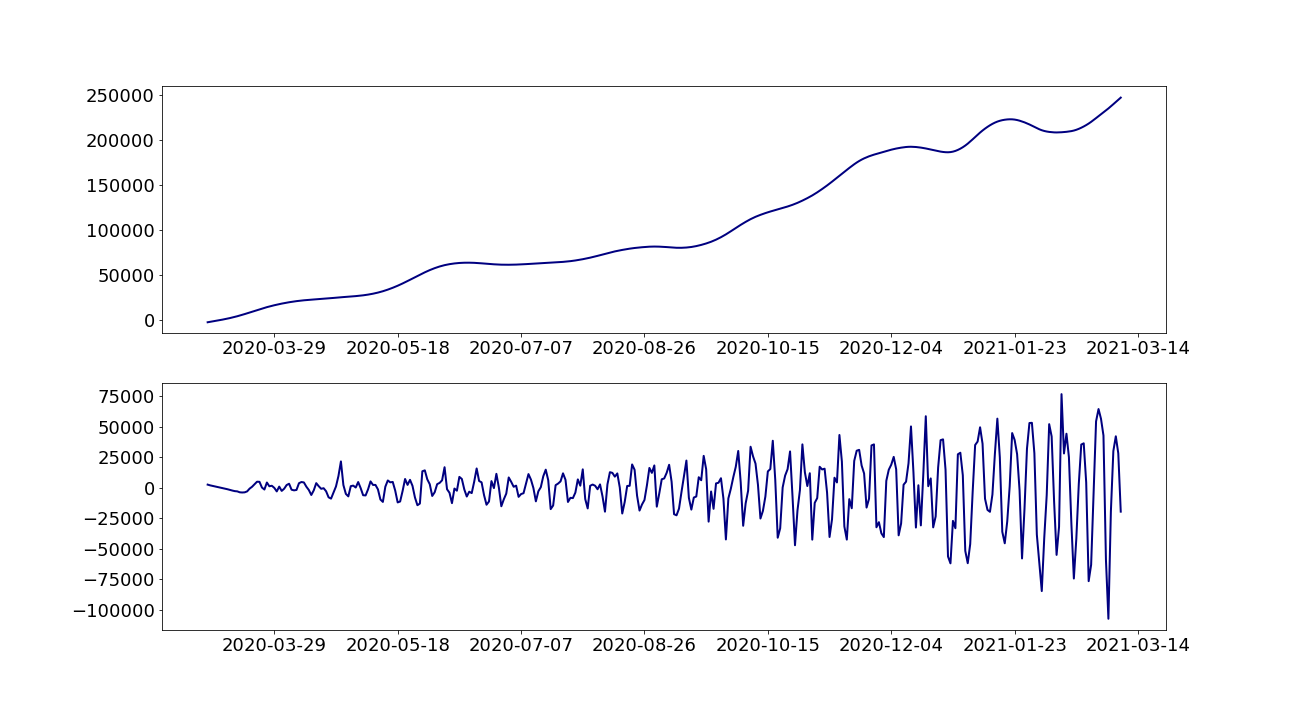}
    \caption{Graphs demonstrating the expansion results of the {\fontfamily{cmss}\selectfont New Tests}: $T(n)$ (top) and $S(n)+N(n)$ (bottom).}
    \label{components}
\end{figure}

We can consider the COVID-19 spread in New York State in more detail, using the epidemic time-series data. Any time series can be decomposed into the following three elements: 

\begin{equation}\label{decom}
    X(n) = T(n) + S(n) + N(n).
\end{equation}
Here, $X(n)$ is a time-series value,  $T(n)$ is the value of the underlying time-series component; $S(n)$ is the value of a seasonality component, $N(n)$ is the value of a noise component for the  $n$-th day. When analyzing a time series, we found it was most useful to analyze its trend, since this determines an indicator’s behavior in time. For expansion by the expression \ref{decom}, the Hodrick-Prescott filter was applied  \cite{smooth} (see result in  Fig.~\ref{components}).

\begin{figure}[H]
    \centering
    \includegraphics[width=\textwidth]{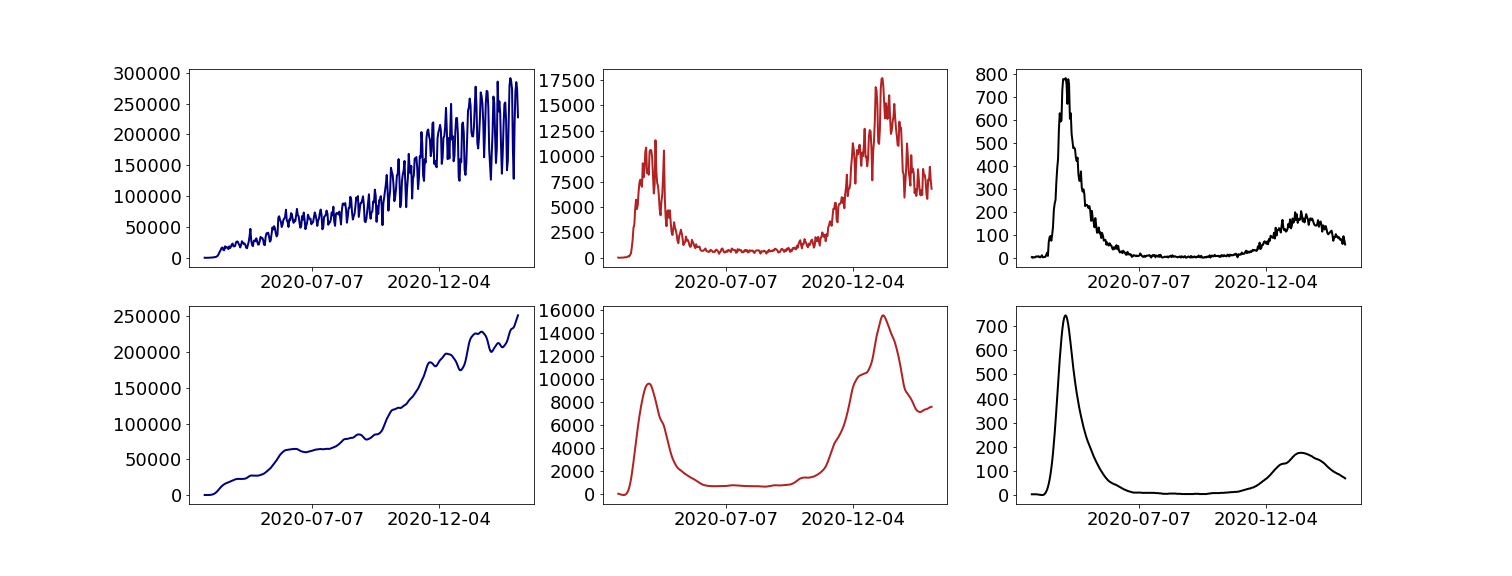}
    \caption{Graphs of COVID-19 spread in New York State (top line): the number of  {\fontfamily{cmss}\selectfont New Tests} (left), {\fontfamily{cmss}\selectfont New Diagnoses} (middle) and {\fontfamily{cmss}\selectfont New Deaths} (right) and their smoothed curves (bottom line).}
    \label{main_graphs_ny}
\end{figure}

The non-smoothed graphs in Fig.~\ref{main_graphs_ny} demonstrate widely dispersed points of statistics, so these data were smoothed before using them in the model because only the main trends of the curves were necessary for reaching an appropriate result. 
The graphs demonstrate certain periodicity that is known to be time series seasonality, which is clearly traced during summer in 
 {\fontfamily{cmss}\selectfont New Diagnoses}.

Now, let us consider the average fraction of tests for each day of the week that is calculated as:

\begin{equation*}
    w_i = \sum\limits_{j=0}^{N} \frac{X(i+m \cdot j)}{S(i+m\cdot j)}.
\end{equation*}
Here, $i$ is a week-day number, $N$ is the number of full weeks within a considered time series, $S(i+m\cdot j)$ is the cumulative sum of the {\fontfamily{cmss}\selectfont New Tests} performed within a week corresponding to index $j$, $m=7$. The results can be seen in Table~\ref{tab:test_data_season_ny}. It is apparent that the number of tests tends to its minimum for Mondays and Sundays and reaches its maximum on Fridays and Saturdays, which exemplifies the seasonality (i.e., periodicity) of the considered time series.

\begin{table}[H]
\begin{center}
\caption{Average fraction from the number of tests and its daily distribution for a whole statistics-gathering period (NY State).}
\label{tab:test_data_season_ny}
\begin{tabular}{|p{4cm}|p{3cm}| } 
\hline
Days of the week & Average fraction\\ 
\hline
Monday & 0.113654 \\
Tuesday & 0.127697\\
Wednesday & 0.134512\\
Thursday & 0.155779\\
Friday & 0.162017\\
Saturday & 0.163702\\
Sunday  & 0.142640\\
\hline
\end{tabular}
\end{center}
\end{table}

COVID-19 pandemic data analysis has shown that there is a dependence between {\fontfamily{cmss}\selectfont New Tests} and {\fontfamily{cmss}\selectfont New Diagnoses} datasets for different regions. In such cases, the behavior of the second indicator is partly determined by that of the first. As a matter of fact, the last day of the week does not mean an abrupt increase in the number of the infected only because it is Friday. However, more people get tested on Friday, so the number of positive tests may increase. For that reason, analyzing the links between the number of infected people and that of new tests becomes crucial. For a proper understanding of the situation, it is not {\fontfamily{cmss}\selectfont New Diagnoses}, but their fraction from the number of {\fontfamily{cmss}\selectfont New Tests} that we want to know. After all, if one tested every person in a region, they would immediately indicate every one infected and their indicator would increase abruptly while their percentage would remain the same. Based on this fact, one can derive a time series that describes the percentage of newly tested people with either positive or potentially positive COVID-19 test.

\begin{figure}[H]
    \centering
    \includegraphics[width=\textwidth]{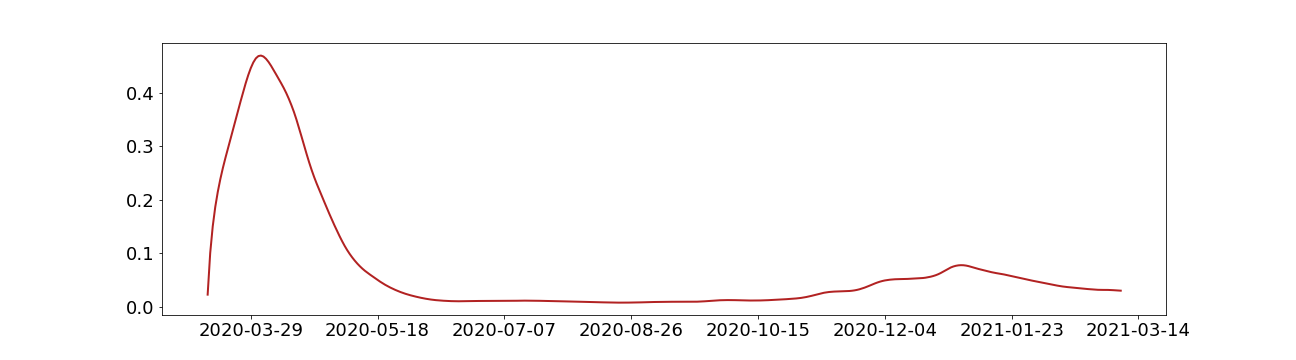}
    \caption{Fraction of {\fontfamily{cmss}\selectfont New Diagnoses} from {\fontfamily{cmss}\selectfont New Tests} in NY State.}
    \label{ratio_ny}
\end{figure}
Figure~\ref{ratio_ny} shows that the ratio reached its peak in April 2020 and later started to reduce due to an increased number of tests. In autumn, the second pandemic wave began, so the number of daily confirmed cases increased (see Fig.~\ref{main_graphs_ny}), but the ratio between {\fontfamily{cmss}\selectfont New Diagnoses} and {\fontfamily{cmss}\selectfont New Tests} remained at a low level. This implies that testing behavior (i.e., the probability of testing with or without COVID-19 symptoms) did not change. The relatively few number of deaths over this period was due to the epidemic primarily spreading among younger age groups (see  Fig.~\ref{main_graphs_ny}).

\begin{figure}[H]
    \centering
    \includegraphics[width=\textwidth]{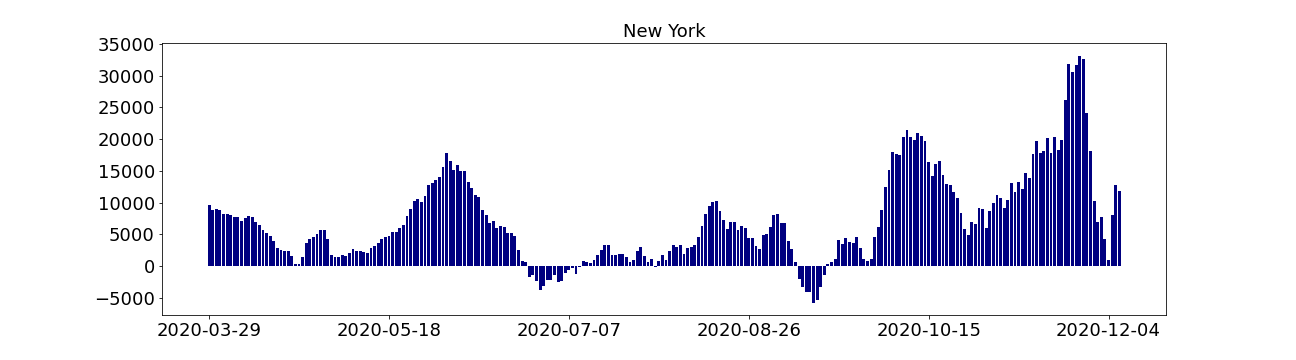}
    \caption{MACD-indicator of the  {\fontfamily{cmss}\selectfont New Tests} in NY State.}
    \label{MACD_ny}
\end{figure}

To confirm change in testing during autumn and winter, the MACD indicator \cite{Bartolucci2018} was used, whose histogram tracks a function’s rise and fall (in our case it was the rise of {\fontfamily{cmss}\selectfont New Tests}). The graph demonstrates that the MACD indicator increased significantly over this period (see Fig.~\ref{MACD_ny}).

\subsubsection{Seasonality}
The above-mentioned term “seasonality” refers to the periodic fluctuations observed in time series. In other words, if one takes a time-series space and overlaps it against a neighboring space of the same size, their profiles are going to coincide (or differ by a constant). The peak absolute values will correlate with the same time points calculated from the beginning of the space. As for our COVID-19 data, their seasonality can be traced in the {\fontfamily{cmss}\selectfont New Tests} indicator, which has its logic since collecting this statistics involves as a human factor as the features of the healthcare system. The seasonality of a time series makes it possible to forecast the series behavior relative to some average value. 

Seasonality is commonly determined with an autocorrelation function (ACF). For a discrete process $X(1), X(2), \ldots, X(n)$ its formula is written as: 
\begin{equation*}
    R(n) = \frac{1}{(M-n) \cdot \sigma^2}\sum\limits_{t=1}^{M-n} (X(t) - \mu)(X(t+n) - \mu).
\end{equation*}
Here $\sigma$ is the standard deviation of a discrete process  $X$, $\mu$ is its average value,  $M$ and $n$ are positive integers. Besides, to analyze time-series seasonality, a partial correlation (PACF) was applied that removed the linear dependence between shifted time series.

\begin{figure}[H]
    \centering
    \includegraphics[width=\textwidth]{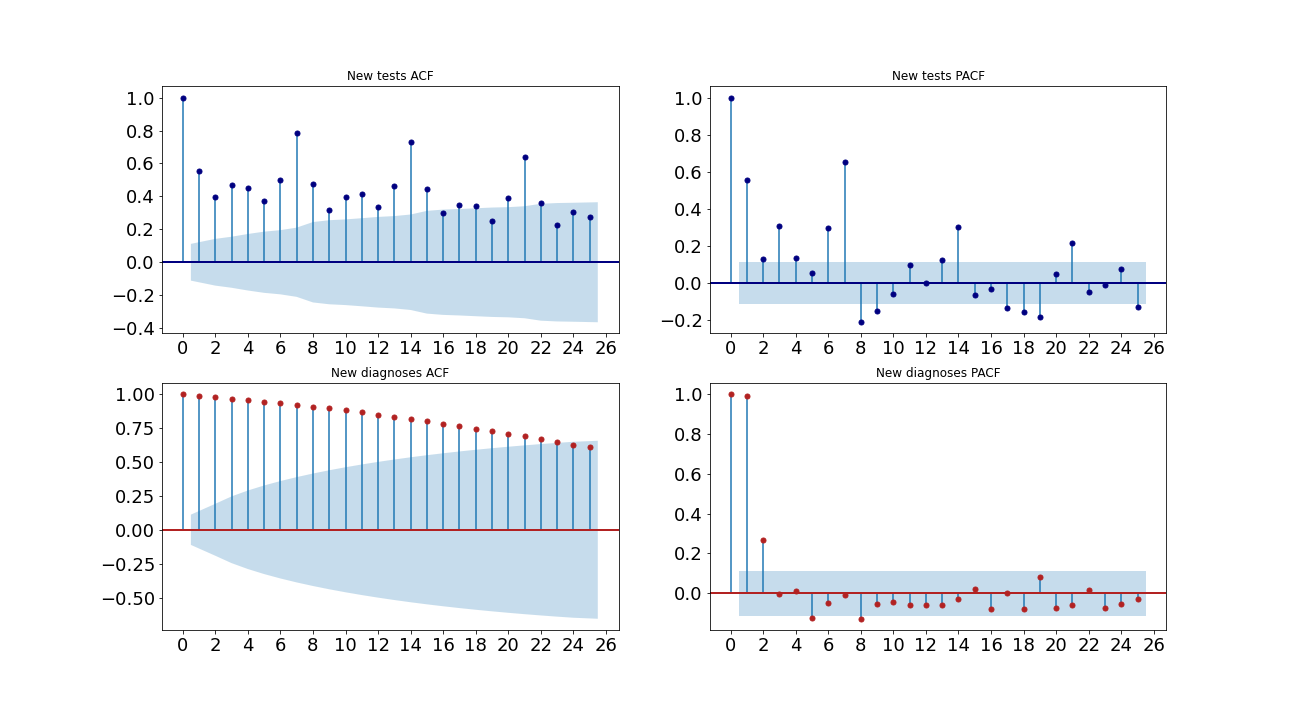}
    \caption{Graphical results of applying ACF (left) and PACF (right) to {\fontfamily{cmss}\selectfont New Tests} (top) and {\fontfamily{cmss}\selectfont New Diagnoses} (bottom) in NY State.}
    \label{autocorr_ny}
\end{figure}

The results confirmed {\fontfamily{cmss}\selectfont New Tests} really had weekly seasonality (7, 14 and 21st days), while the {\fontfamily{cmss}\selectfont New Diagnoses} did not possess this property (see Fig.~\ref{autocorr_ny}).

\subsubsection{New Tests/New Diagnoses interrelation}
We have also considered a percentage change for any current moment in relation to the same moment a week before:
\begin{equation*}
    pc(n) = {X_{smoothed}(n) \over X_{smoothed}(n-7)} - 1.
\end{equation*}
This statistic demonstrates by how much an average indicator has changed in fractions compared to a previous week. If compared for weekdays, the results become smoother and easier to interpret. 
\begin{figure}[H]
    \centering
    \includegraphics[width=\textwidth]{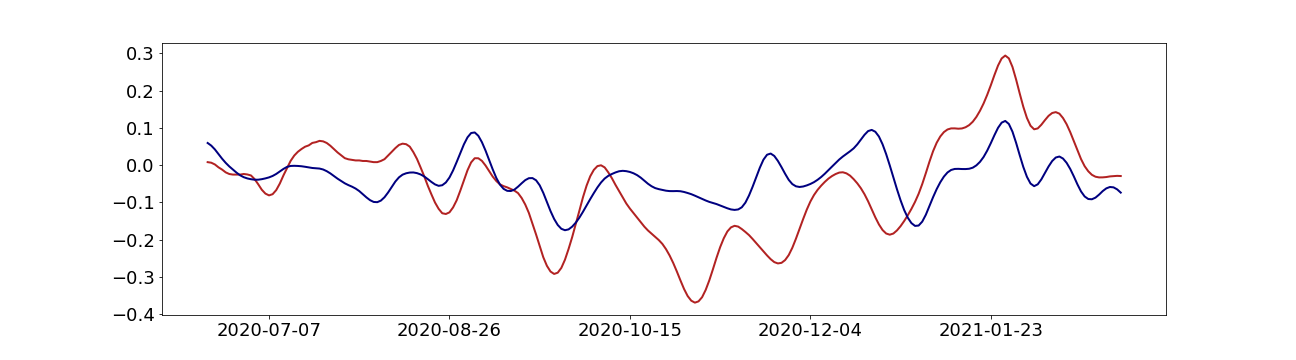}
    \caption{{\fontfamily{cmss}\selectfont New Diagnoses} percent change (red line) and  {\fontfamily{cmss}\selectfont New Tests} percent change (blue line) time-series graphs in NY State.}
    \label{pc_ny}
\end{figure}
The two time series in Fig.~\ref{pc_ny} have spaces where their trajectories almost match. In other words, the percent change of one indicator differs from that of the other one by a constant. In terms of tested/infected ratio, such spaces confirm that within this period, the number of infected people grew owing to the increased number of tests and not to a worsening pandemic situation in the region. 

\begin{figure}[H]
    \centering
    \includegraphics[width=\textwidth]{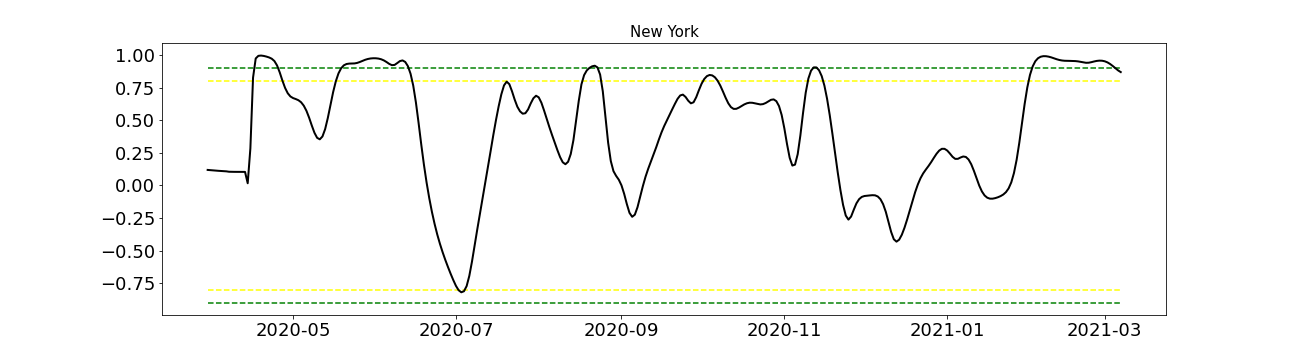}
    \caption{{\fontfamily{cmss}\selectfont New Diagnoses} percent change and {\fontfamily{cmss}\selectfont New Tests} perecent change correlations within a 28-day backward window with $\mid0.8\mid$ and $\mid0.9\mid$ correlation's threshold lines for NY state.}
    \label{corr_ny}
\end{figure}

For more specific analysis, a window of 28 days was correlated with the previous 28 days for every day in the window  (see Fig.~\ref{corr_ny}).
In statistics, a significant linear dependence occurs when the absolute correlation value exceeds 0.8. If it exceeds 0.9, it means there is a strong dependence between two indicators. In NY State, with some deviations,  {\fontfamily{cmss}\selectfont New Diagnoses} strongly depended on  {\fontfamily{cmss}\selectfont New Tests}. The deviations mean the pandemic develops following its own scenario. For instance, in the end of May and the beginning of June, the correlation coefficient exceeded 0.8 within a 28-day backward window when the first pandemic wave in the region was considered defeated and reached the so-called plateau that lasted till October when the second correlation took place and the second pandemic wave began, which means the second wave might have been triggered by the abrupt increase of {\fontfamily{cmss}\selectfont New Tests}. However, the correlation coefficient got back to almost zero by December, indicating that the rise and fall of {\fontfamily{cmss}\selectfont New Diagnoses} did not correlate with  {\fontfamily{cmss}\selectfont New Tests}, so the pandemic in the region at that time spread or reduced (depending on graph direction). In February and March, the correlation became very close to one, meaning the number of infected people depended only on a {\fontfamily{cmss}\selectfont New Tests}.

\subsection{Novosibirsk region}
\begin{figure}[H]
    \includegraphics[width=12cm]{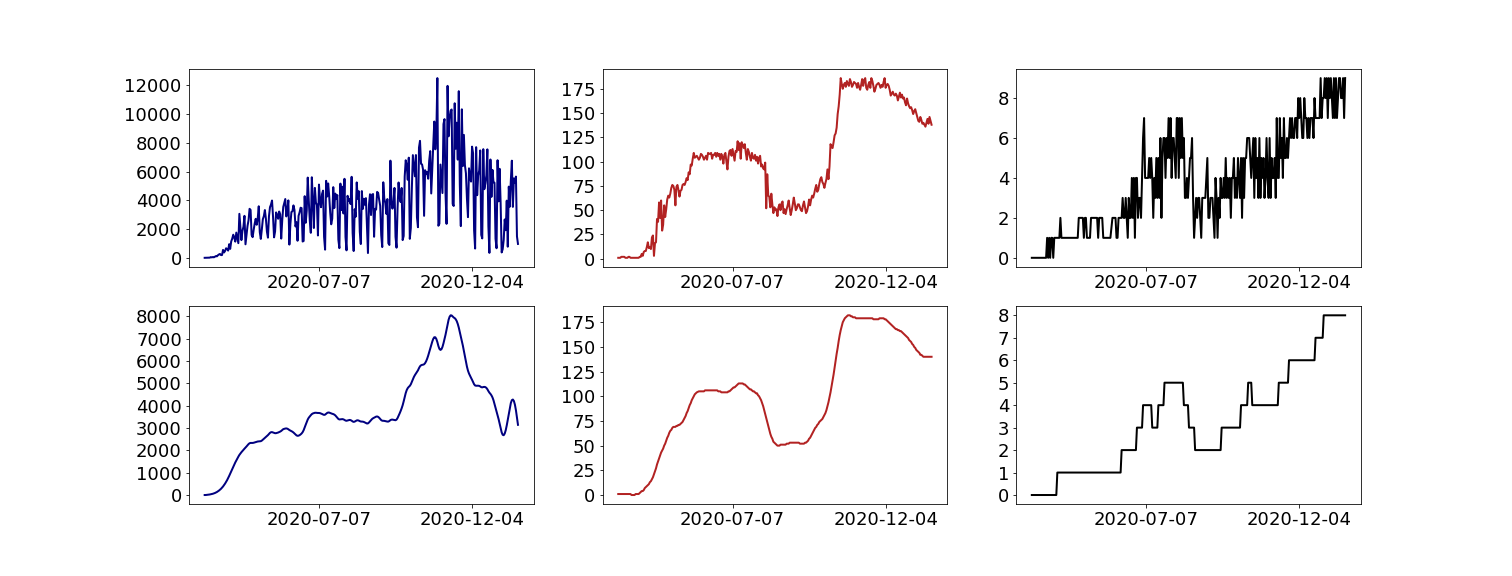}
    \caption{COVID-19 spread graphs for the Novosibirsk region (top line): {\fontfamily{cmss}\selectfont New Tests} (left), {\fontfamily{cmss}\selectfont New Diagnoses} (middle) and {\fontfamily{cmss}\selectfont New Deaths} (right) and their smoothed curves (bottom line).}
    \label{main_graphs_nsk}
\end{figure}

\begin{figure}[H]
    \includegraphics[width=12cm]{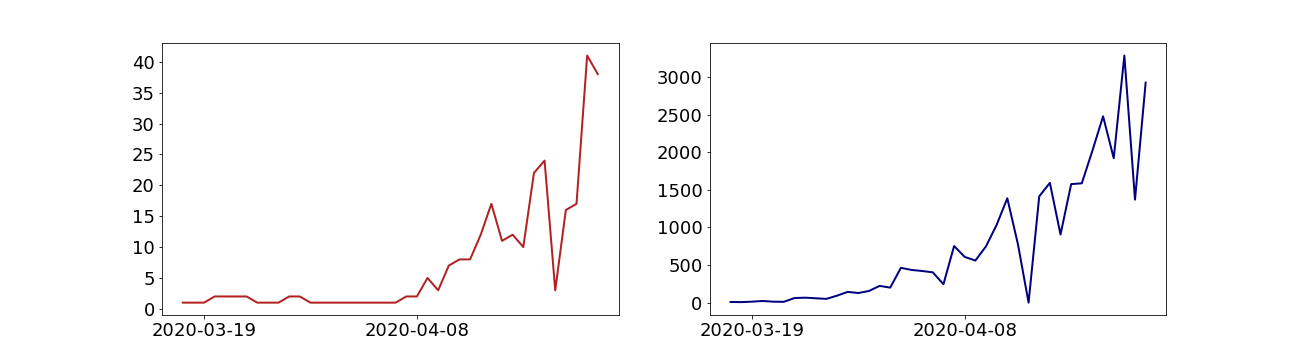}
    \caption{Backward extrapolation results for the {\fontfamily{cmss}\selectfont New Tests} (right) compared with the actual values of {\fontfamily{cmss}\selectfont New Diagnoses} time series in the Novosibirsk region (left).}
    \label{tests_extr}
\end{figure}

The data (see Fig.~\ref{main_graphs_nsk}) for our research were taken from the RBC website that published the daily statistics of COVID-19 spread in the Novosibirsk region (the used data are given in Table \ref{tab:data_statistics}). In March and May of 2020, these data had a significant number of gaps. For instance, the statistic on the number of performed tests was available only starting April 25, 2020, as the other statistics were accumulated starting March 17, 2020. For that reason, the {\fontfamily{cmss}\selectfont New Tests} was extrapolated in relation to the abrupt increase in the number of {\fontfamily{cmss}\selectfont New Diagnoses} starting April 6, 2020. For filling the gaps, the following backward extrapolation algorithm was iterated:

\begin{equation*}
    X(n) = \frac{1}{m}\sum\limits_{j=1}^{m} X(i+j)\cdot C^n + E(n), \quad i=L,\ldots ,0.
\end{equation*}
Here, $X(n)$ is the value of a time series of the $n$-th day; $m = 7$, $C = 1.03$, $L$ is the number of gaps to be extrapolated; $E(n)$ is a random value of normal distribution, i.e. $E(n)\in \mathcal{N}(0, {X(n) \over 3})$. The extrapolation result for {\fontfamily{cmss}\selectfont New Tests} can be seen in Fig.~\ref{tests_extr}.

\subsubsection{Seasonality}
\begin{figure}[H]
    \centering
    \includegraphics[width=\textwidth]{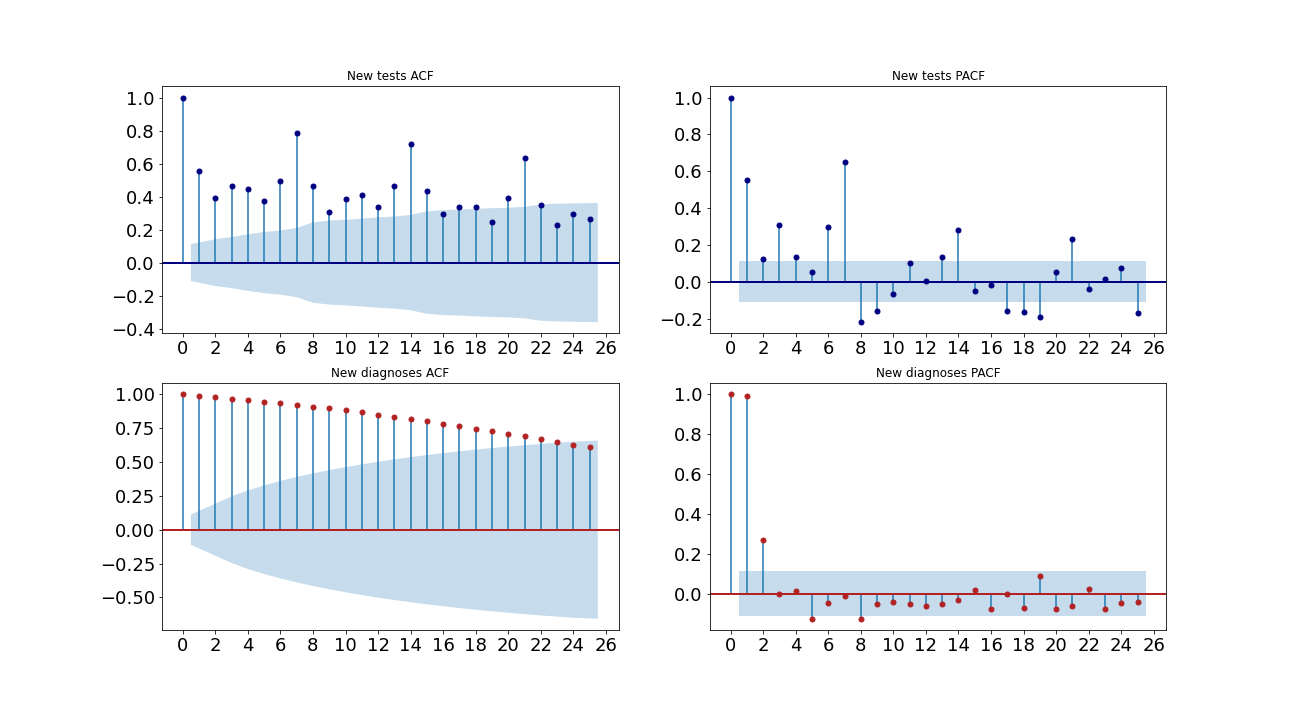}
    \caption{Graphical results of applying ACF (left) PACF (right) to the {\fontfamily{cmss}\selectfont New Tests} (top) and {\fontfamily{cmss}\selectfont New Diagnoses} (bottom) in the Novosibirsk region.}
    \label{autocorr_nsk}
\end{figure}

Here, the approach similar to that applied for the NY State data was implemented to calculate the ACF and PACF for several time series (see Fig.~\ref{autocorr_nsk}), together with the fraction of tests performed for every day of the week (see Table~\ref{tab:test_data_season_nsk}).

\begin{table}[H]
\begin{center}
\caption{Average fraction from the number of tests and its daily distribution for a whole statistics-gathering period in the Novosibirsk region.}
\label{tab:test_data_season_nsk}
\begin{tabular}{|p{4cm}|p{3cm}| } 
\hline
Week days & Average fraction\\ 
\hline
Monday & 0.064055 \\
Tuesday & 0.176845\\
Wednesday & 0.148333\\
Thursday & 0.161924\\
Friday & 0.162173\\
Saturday & 0.183951\\
Sunday  & 0.075146\\
\hline
\end{tabular}
\end{center}
\end{table}

As can be seen from Fig.~\ref{autocorr_nsk}, both functions confirmed the seasonality of {\fontfamily{cmss}\selectfont New Tests} with the duration of the season to be one week. As it has been discussed earlier, for this indicator this is absolutely normal. Table~\ref{tab:test_data_season_nsk} demonstrates that the lowest fraction of tests was registered on Sunday and Monday, and the highest one – on Tuesday and Saturday. 

\subsubsection{Stationarity}\label{sec_NSK_stationarity}
A time series is regarded as stationary (or weakly stationary) if its expected value and variance are time-independent, and its autocorrelation function depends only on a difference of neighboring values. Let us consider two time series for the {\fontfamily{cmss}\selectfont New Diagnoses}: from May 29 to July 1 and from October 24 to December 7 (see Fig.~\ref{main_graphs_nsk}). Having applied the augmented Dickey-Fuller test, whose zero hypothesis is non-stationary time series (Unit root test), we find that this hypothesis can be rejected for both time series intervals, since their p-values are 0.0 and 0.001, respectively. If time series are stationary, the considered data sets can probably be interpreted as a normal noise. Table~\ref{tab:stationary} demonstrates the results of two statistical tests to verify the normality of the considered time series. 
\begin{table}[H]
\begin{center}
\caption{Results of the two statistical tests, whose zero hypothesis is the normal distribution of two stationary parts of {\fontfamily{cmss}\selectfont New Diagnoses}.}
\label{tab:stationary}
\begin{tabular}{|p{3.6cm}|p{3.8cm}|p{3.8cm}| } 
\hline
Time interval & Kolmogorov-Smirnov test ($p$-value) & D'Agostino-Pearson test ($p$-value) \\ 
\hline
29.05.2020 - 1.07.2020 & 0.3842 & 0.1827\\
24.10.2020 - 7.12.2020 & 0.573 & 0.8193\\
\hline
\end{tabular}
\end{center}
\end{table}
According to the test results, the null hypothesis of the normality of the two areas is not rejected. However, due to the small number of observations, it is equally impossible to assert that the values in the selected areas are normally distributed.

\subsubsection{New Tests/New Diagnoses interrelation}
\begin{figure}[H]
    \centering
    \includegraphics[width=\textwidth]{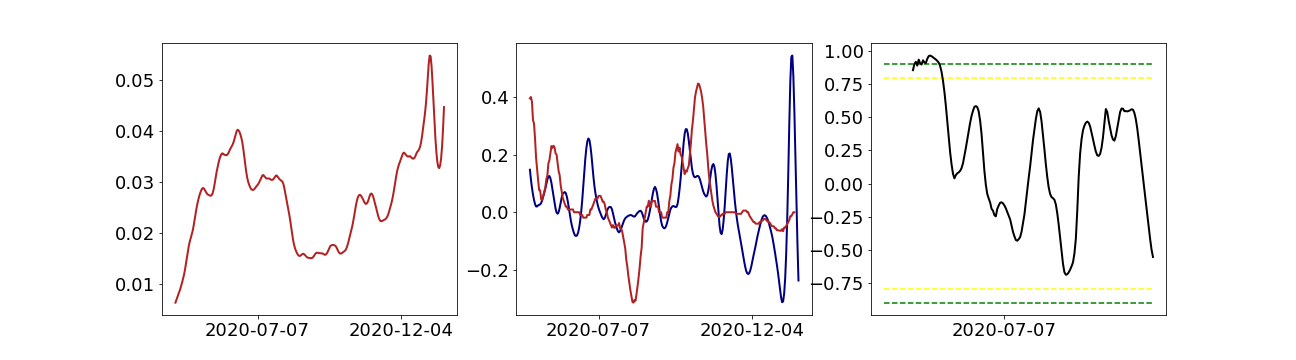}
    \caption{{\fontfamily{cmss}\selectfont New Diagnoses} to {\fontfamily{cmss}\selectfont New Tests} ratio (left); {\fontfamily{cmss}\selectfont New Diagnoses} (red) percent change and {\fontfamily{cmss}\selectfont New Tests} (blue) percent change (center) and their correlations within a 28-day backward window (right) with $\mid 0.8 \mid$ and $\mid 0.9 \mid$ correlation’s threshold lines for the Novosibirsk region.}
    \label{Novosib_stats}
\end{figure}

Calculated {\fontfamily{cmss}\selectfont New Diagnoses} and {\fontfamily{cmss}\selectfont New Tests} percent change correlations (see Fig.~\ref{Novosib_stats}) showed there was not even a single window of 28 days with a significant {\fontfamily{cmss}\selectfont New Tests}/{\fontfamily{cmss}\selectfont New Diagnoses} linear dependence. In the middle of June, a negative correlation peak was observed confirming the “more tests less positive cases” tendency lasting from the end of June to the beginning of July. Another interesting fact was the {\fontfamily{cmss}\selectfont New Diagnoses}, percentage that plummeted in August while {\fontfamily{cmss}\selectfont New Tests} remained at the same level.

\subsection{United Kingdom} \label{sect: UK}
To analyze COVID-19 spread in the UK, the data accumulated on an official government web portal were used \cite{uk_gov}.

\begin{figure}[H]
    \centering
    \includegraphics[width=\textwidth]{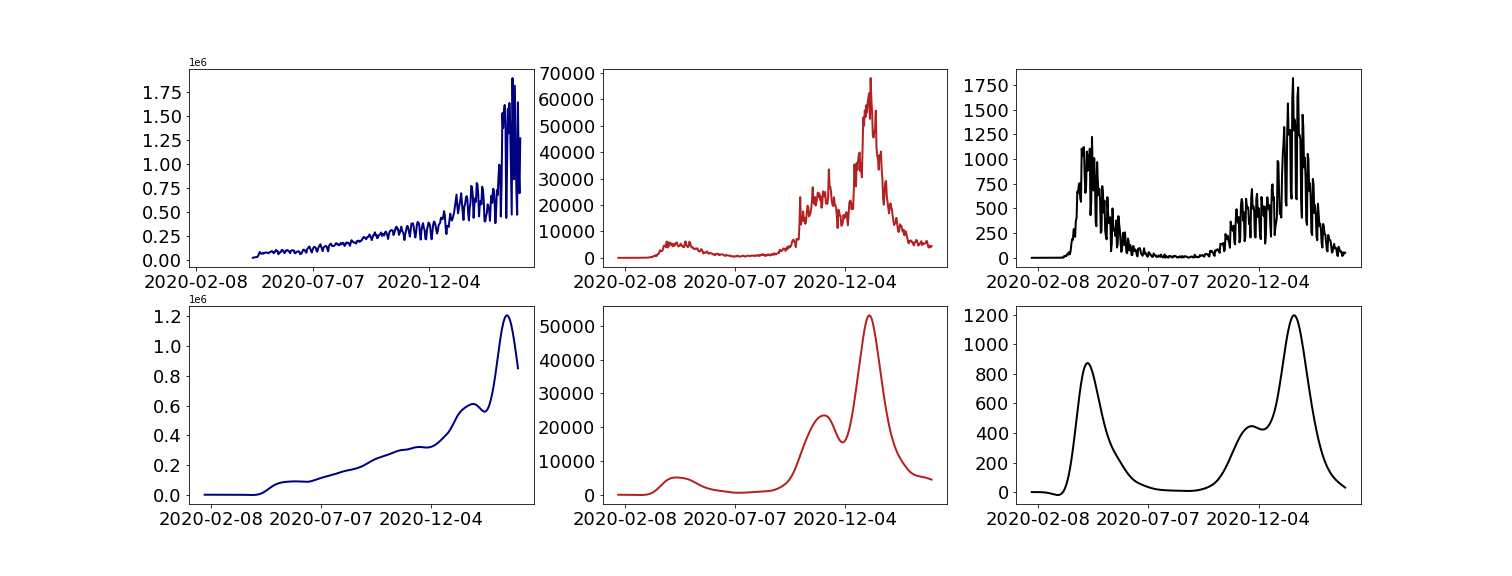}
    \caption{ COVID-19 spread graphs in the UK (top line): {\fontfamily{cmss}\selectfont New Tests} (left), {\fontfamily{cmss}\selectfont New Diagnoses} (middle),  {\fontfamily{cmss}\selectfont New Deaths} (right) and their smoothed curves (bottom line).}
    \label{main_graphs_uk}
\end{figure}

Fig.~\ref{main_graphs_uk} shows two infection waves: spring (from 2020.03.09 to 2020.06.01) and fall (from 2020.09.01 to 2021.03.01) ones. The standard deviation of the time series for the second wave was much bigger than for the first one. For that reason, the standard deviation values were considered for three independent time series of  {\fontfamily{cmss}\selectfont New Diagnoses} and {\fontfamily{cmss}\selectfont New Tests} (see Table~\ref{tab:data_statistics_uk}). The data for NY State (NY) and the Novosibirsk region (NSK) are put in the table for comparison.

\begin{table}[H]
\begin{center}
\caption{Standard deviation values for particular {\fontfamily{cmss}\selectfont New Diagnoses} and {\fontfamily{cmss}\selectfont New Tests} time series in the UK, NY State and the Novosibirsk region.}
\label{tab:data_statistics_uk}
\begin{tabular}{|p{4cm}|p{1.2cm}|p{1.2cm}|p{1.2cm}|p{1.2cm}|} 
\hline
Indicator & First & Second  & Third & Fourth\\ 
\hline
(UK) {\fontfamily{cmss}\selectfont New Tests std} & 10626.7 & 16613.8 & 77382.5 & ...\\
(UK) {\fontfamily{cmss}\selectfont New Diagnoses std} & 455.8 & 201.1 & 5727.2 & ... \\
\hline
(NY) {\fontfamily{cmss}\selectfont New Tests std} & 5636.6 & 9543.9 & 31127 & ... \\
(NY) {\fontfamily{cmss}\selectfont New Diagnoses std} & 1137 & 114.2 & 1144 & ... \\
\hline
(NSK) {\fontfamily{cmss}\selectfont New Tests std} & 765.2 & 1271.2 & 1596 & 2492 \\
(NSK) {\fontfamily{cmss}\selectfont New Diagnoses std} & 7 & 4.6 & 7 & 6 \\
\hline
\end{tabular}
\end{center}
\end{table}

Unfortunately, the data did not allow us to conclude what was the exact reason for the second wave’s higher variance. It could have been the week variance of  {\fontfamily{cmss}\selectfont New Tests} or something else since, in the UK, they started to register the number of tests only on 2020.04.21. All we know is the standard deviation of  {\fontfamily{cmss}\selectfont New Diagnoses} increased with time as well as the variance of  {\fontfamily{cmss}\selectfont New Tests}. However, both in NY State and the Novosibirsk region, despite the growing variance of  {\fontfamily{cmss}\selectfont New Tests}, the standard deviation value of  {\fontfamily{cmss}\selectfont New Diagnoses} remained at the same level.

\subsubsection{Seasonality}
\begin{figure}[H]
    \centering
    \includegraphics[width=\textwidth]{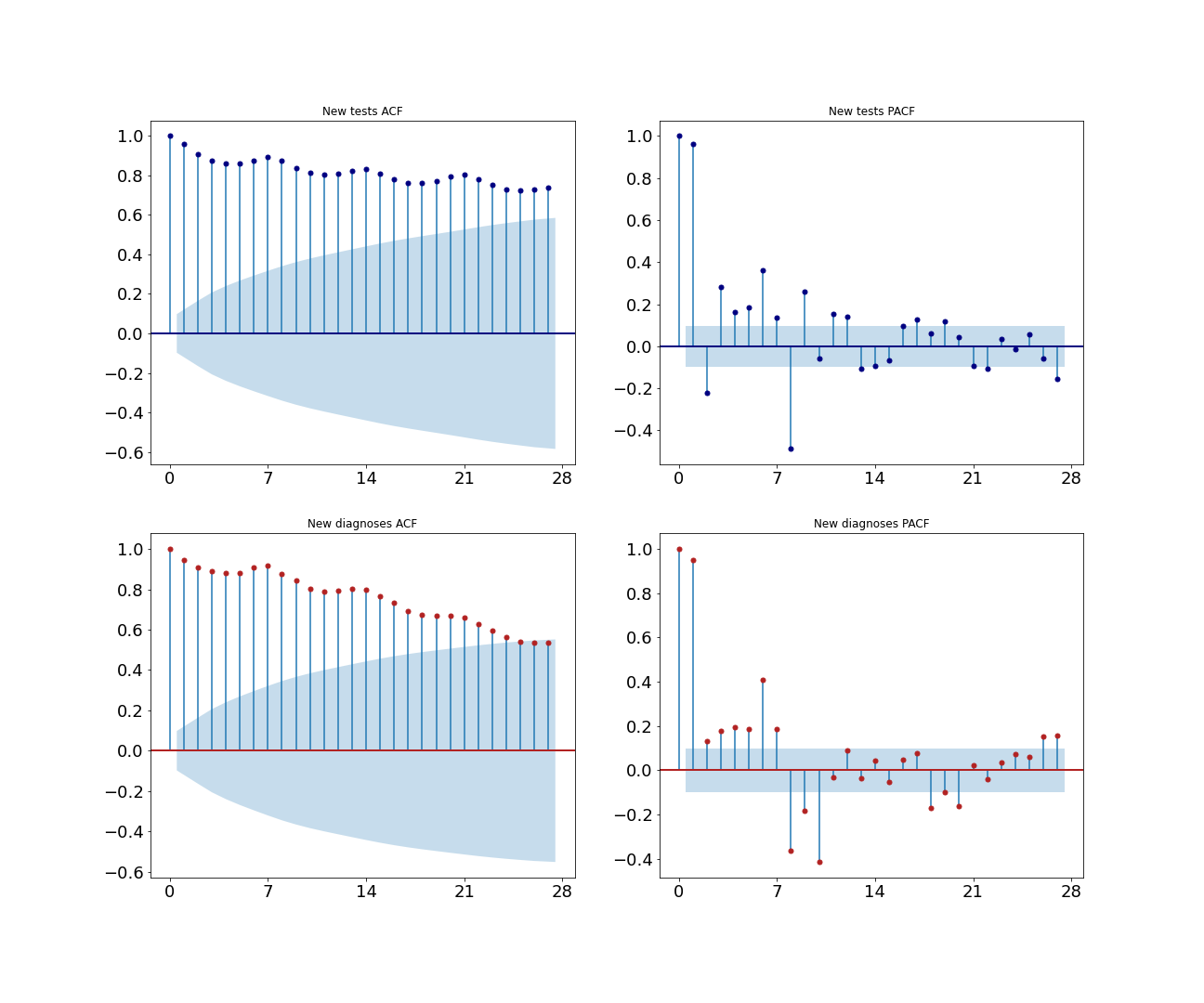}
    \caption{Graphical results of applying ACF (left) and PACF (right) to {\fontfamily{cmss}\selectfont New Tests} (top) and {\fontfamily{cmss}\selectfont New Diagnoses} (bottom) for the UK.}
    \label{autocorr_uk}
\end{figure}

The ACF and PACF applied to the {\fontfamily{cmss}\selectfont New Diagnoses} and {\fontfamily{cmss}\selectfont New Tests} time series and the fractions of tests calculated for every day of the week demonstrated the weekly seasonality of {\fontfamily{cmss}\selectfont New Tests} confirmed by ACF (see Fig.~\ref{autocorr_uk} and Table~\ref{tab:test_data_season_uk}). However, no significant PACF delays were observed. Comparing the results obtained with those from NY State and the Novosibirsk region demonstrated the absence of linear dependence did not affect the time series’s seasonality. However, a different situation was observed in the UK.

\begin{table}[H]
\begin{center}
\caption{Average fraction from the number of tests and its daily distribution for a whole statistics-gathering period in the UK.}
\label{tab:test_data_season_uk}
\begin{tabular}{|p{4cm}|p{3cm}| } 
\hline
Days of the week & Average fraction\\ 
\hline
Monday & 0.125222 \\
Tuesday & 0.135905 \\
Wednesday & 0.159167 \\
Thursday & 0.166208 \\
Friday & 0.157837 \\
Saturday & 0.131992 \\
Sunday  & 0.123668 \\
\hline
\end{tabular}
\end{center}
\end{table}

In the next section, an approach to forecasting the {\fontfamily{cmss}\selectfont New Tests} time series will be demonstrated. The approach relies upon several techniques, including SARIMA, an algorithm requiring a significant correlation for the seasonal delays (7, 14, 21, etc.) of a single parameter.

\subsection{Forecasting New Tests}
To draw forecast curves while modeling, one has to predict a number of certain statistical data sets that are used as input parameters. In our model, such a data set was the {\fontfamily{cmss}\selectfont New Tests} performed in the region since this indicator did not depend on the others and had the highest value of seasonality that determined the seasonality of the other indicators.

The first model considered for forecast purposes was SARIMA~\cite{sarima}, a modification of ARIMA~\cite{arima} (AutoRegressive Integrated Moving Average), which is able to support time-series data sets with a seasonal component. ARIMA is an extension of the ARMA models for non-stationary time series that can be converted into stationary ones through differencing of a certain order from an initial time series (so-called integrated or differential-stationary time series). 

For a non-stationary time series $X(n)$, the ARIMA$(p,d,q)$ model is written as:
\begin{equation*}
    \triangle^d X(n) = c + \sum_{i=1}^p a_i \triangle^d X(n-i) + \sum_{j=1}^p b_j \epsilon(n-j) + \epsilon(n).
\end{equation*}
Here, $\epsilon(n)$ is a stationary time series of white noise; $c, a_i, b_j$ -- model parameters; $\triangle^d$ -- difference operator of a time series of order $d$ (sequential taking of $d$ times of first-order differences: firstly, from a time series, then from obtained first-order differences, then from second-order differences and so on). 

Parameters $p$, $d$, $q$ and $P$, $D$, $Q$ are determined through cross-validating minimization of the AIC-metrics. Here,  $p$ is the number of the last non-seasonal delay with significant PACF;  $P$ is the number of the last seasonal delay with significant PACF; $q$ is the number of the last non-seasonal delay with significant ACF;  $Q$ is the number of the last seasonal delay with significant ACF; $d$ is a differentiation order with 1-day delay, $D$ is a seasonal differentiation order with 7-day delay.  The ARIMA model is based on the following idea: one builds a model from the differences of a value for $d$ sequential periods to obtain a stationary process. In a general case, the difference order is limited to $d=2$, since taking the second-order differences allows converting almost any non-stationary data series into stationary ones. 

As for the time-series forecast, the following algorithm was applied: 
\begin{enumerate} 
  \item Applying the Box-Cox transformation \cite{boxcox} to reduce the variance. 
  \item Calculating the first-order seasonal difference (7-day shift).
  \item Calculating the difference (1-day shift) of the series obtained at step 2.
  \item Applying the Dickey-Fuller criterion to verify the stationarity of the series obtained at step 3. 
  \item Passing the hyperparameters corresponding to the actions performed and selecting other ones using the minimized Akaike information criterion. As data, the series from step 1 is passed.
  \item The model with adjusted hyperparameters was used for forecasting. The obtained results were treated using the inverted Box-Cox transform.
\end{enumerate}

The second model considered for the forecast was Holt-Winters \cite{hw}, which provides triple exponential smoothing of a time series. The idea behind the model is to break a time series into such components as a level, trend and seasonality. The last component explains repeating fluctuations around the level and the trend and is characterized by its duration, i.e. a period after which a fluctuation repeats itself. In this approach, each observation has its own component, e.g. if a season's duration is 7 days, there are 7 seasonal components, one for every day of the week. 

To build the model, parameters $\alpha$, $\beta$  and  $\gamma$ had to be determined. $\alpha$ is responsible for smoothing the series around the trend, $\beta$ is for smoothing the trend itself, $\gamma$ is the seasonality component. Minimizing mean absolute error (MAE) through time-series cross-validation, we selected the model’s optimal parameters (iterative three-parameter random search). 

\begin{figure}[H]
    \centering
    \includegraphics[width=\textwidth]{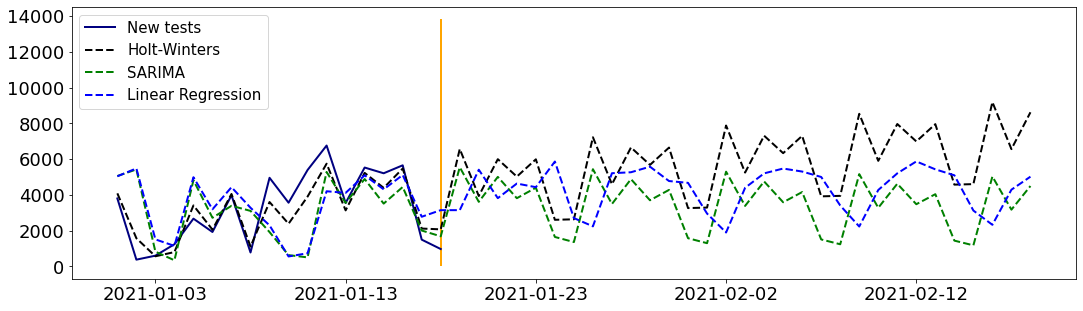}
    \caption{Monthly forecast results for the {\fontfamily{cmss}\selectfont New Tests} time series (dark blue line) for 2021.01.18 obtained using the Holt-Winters (black dashed line), linear regression (blue dashed line) and SARIMA (green dashed line) models for the data from the Novosibirsk region.}
    \label{forecasting}
\end{figure}

To extrapolate the {\fontfamily{cmss}\selectfont New Tests} time series from the Novosibirsk region, in addition to the models described above, a linear regression model was used for comparative analysis. In this respect, for different time intervals, SARIMA’s forecast results were on average not better but more stable than those produced by the other models (see Fig.~\ref{forecasting}).

\begin{table}[H]
\begin{center}
\caption{Mean Absolute Error (MAE) between the {\fontfamily{cmss}\selectfont New Tests} time series and a corresponding prediction (from May 5, 2020 to January 18, 2020).}
\label{tab:sarima_mae}
\begin{tabular}{|p{5cm}|p{5cm}|} 
\hline
Forecast model & MAE\\ 
\hline
Holt-Winters model & 1646.49 \\
SARIMA model & 1454.73\\
Linear Regression model & 1586.34\\
\hline
\end{tabular}
\end{center}
\end{table}

As a result, the SARIMA forecast gave the lowest MAE (see Table~\ref{tab:sarima_mae}). This algorithm also produced the most accurate seasonality and time-series autoregression. For that reason, a decision was made to use SARIMA for forecasting the future values of {\fontfamily{cmss}\selectfont New Tests}. 

%% file: 4agent_model.tex
\section{Agent-based mathematical model}
This section presents the ABM devised to describe COVID-19 spread and formulates a problem to identify (calibrate) the model’s unknown parameters as objective functional minimization. It also presents a scheme for automatic calibration of the parameters for time intervals exemplified by the number of infected people in the Novosibirsk region. 

\subsection{ABM formulation}
Within the framework of this research, stochastic ABMs for New York State, the United Kingdom and the Novosibirsk region were devised. They were built using the Covasim  \cite{covasim}. This package had been utilized to predict the number of infected, dead and hospitalized people in the State of Oregon and become one of the tools to make decisions about whether to relax or escalate COVID-19 containment measures \cite{oregon-2,oregon-1}. This library is written in Python to study non-trivial COVID-19 dynamics. Its general algorithm is as follows: after all necessary parameters and statistical data are uploaded, the package creates an artificial population with account for age distribution. The model’s agent is a person in a particular region. Then, the agents are united into contact networks and the integration loop begins. At every time step (1 day), an agent’s status is updated in relation to its contact network and the containment measures relevant for this interval (self-isolation; closed access to public places; wearing face masks, etc.). The agents can interact with one another in particular networks. Depending on the network’s structure, both full and random connectivity graphs are built for an agent to figure out how the infection spreads. The average number of daily contacts is different for every agent and every network. At any time moment, the agents distributed by their age (bins of 0-9, 10-19, \ldots, 90+) are found in their given state (see fig.~\ref{states}). Our ABM also accounted for tested agents according to the data on COVID-19 tests performed in the Novosibirsk region. More details about the structure of Covasim-based ABMs, their parameters and realization methods can be found in~\cite{covasim_doc}.

Every agent has their set of properties and characteristics that can be divided into 2 groups: constant (belong to each particular agent and do not change while modeling) and time-dependent.

\subsubsection{Time–independent agent characteristics}
\begin{itemize}
    \item Age($t^*$). All the agents are subdivided into age groups of 10 years (0-9 years, 10-19, \ldots, 90$+$. The age distribution depends on the demographic situation in a studied region. 

    \item Social status (determined by an agent’s age  $t^*$). Depending on their age, agents contact one another in contact networks. All agents have contacts in households and public places. Agents of 6-21 years old can also have contacts in educational institutions with agents of their age. Agents of 22-65 years old contact at work (see Fig.~\ref{contact_layers}). Depending on a contact’s structure, the transmission parameter $\beta$ is multiplied by corresponding constant $w_{\beta}$ ($w_{\beta}=3$ for households, 0.6 – for educational institutions, 0.3 – for public places), i.e. the likelihood of virus transfer is different for every network.

    \begin{figure}[h]
        \centering
        \includegraphics[width=12cm]{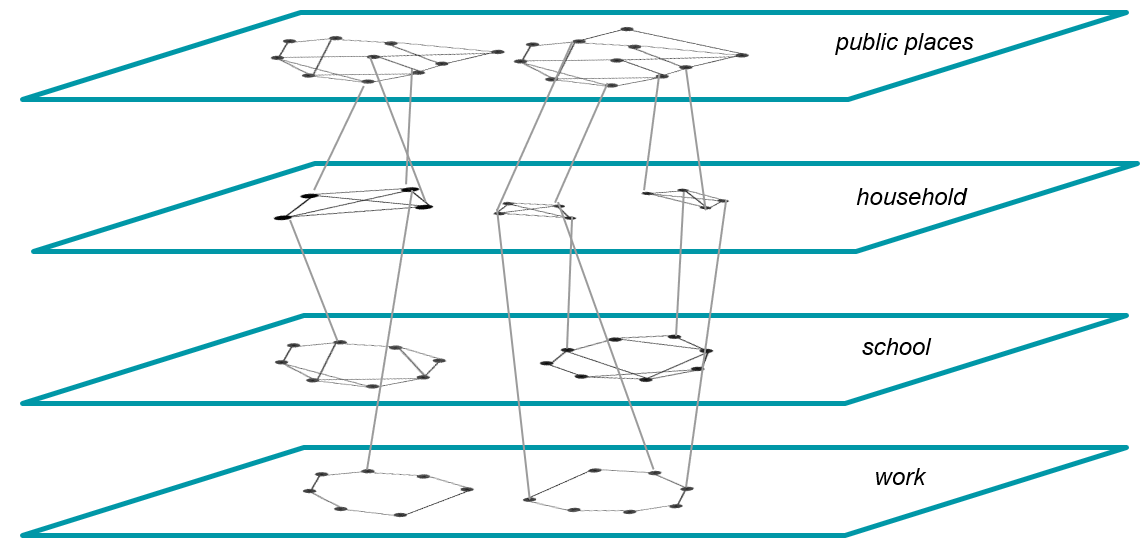}
        \caption{Agent contacts layers and their interactions in the COVID-19 spread ABM.}
        \label{contact_layers}
    \end{figure}

    \item Likelihood of disease progression (determined by an agent’s age  $t^*$). These parameters determine disease progression (see Fig.~\ref{states}). Their description is given in Table~\ref{prob_desease}.
\end{itemize}

\begin{table}[H]
\begin{center}
\caption{Parameters of disease progression probability.}
\label{prob_desease}
\begin{tabular}{|p{4cm}|p{7cm}| } 
\hline
Parameter & Description\\
\hline
$p_{sym}(t^*)$ & Probability of developing symptoms\\ 
$p_{sev}(t^*)$ & Probability of developing severe symptoms (requires hospitalization)  \\
$p_{crit}(t^*)$ & Probability of critical condition
(requires ICU) \\
$p_{death}(t^*)$ &  Probability of death \\

\hline
\end{tabular}
\end{center}
\end{table}

\subsubsection{Time-dependent agent characteristics}
\begin{figure}[ht]
    \centering
    \includegraphics[width=12cm]{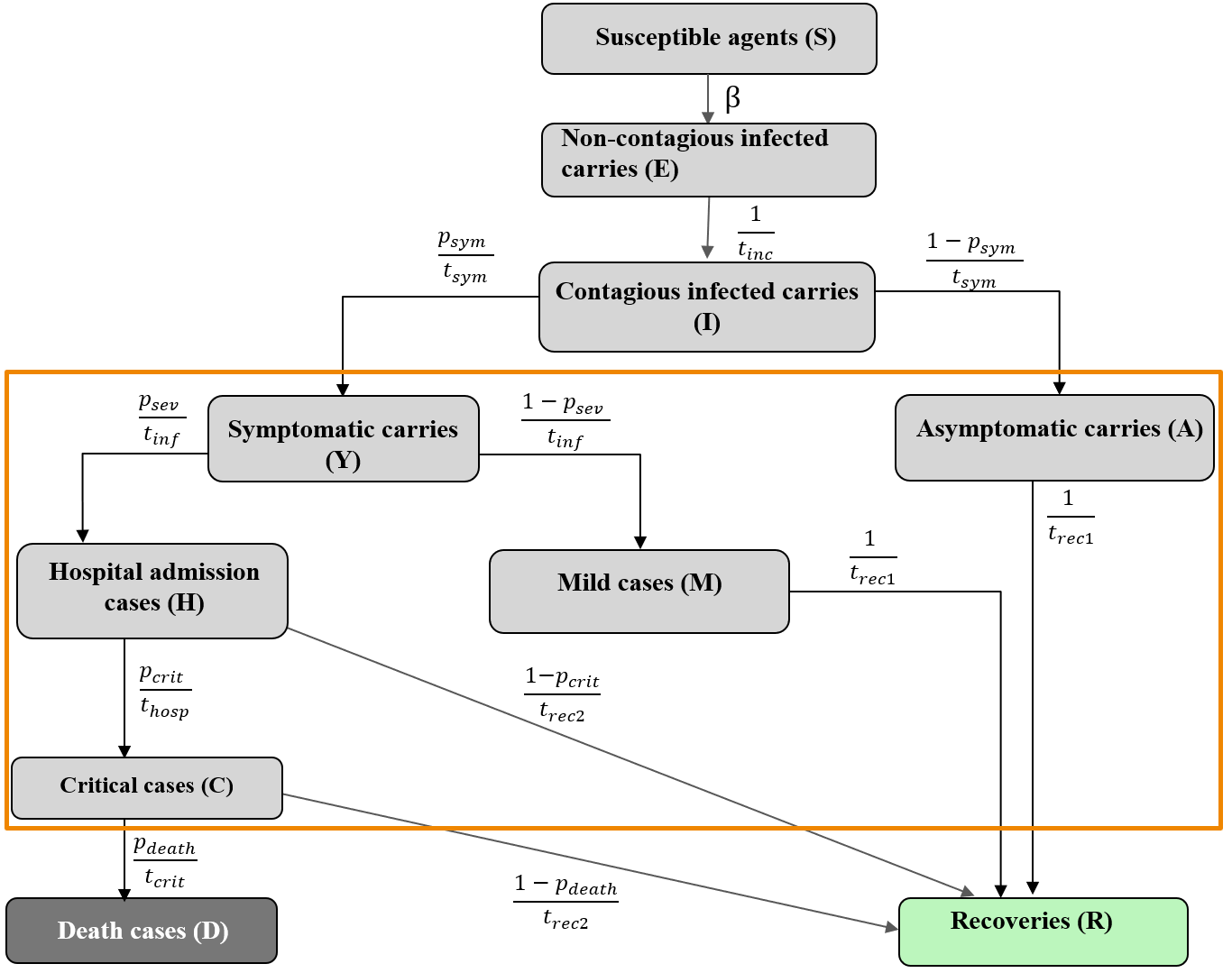}
    \caption{Agent state transition diagram in Covasim that is based on a SEIR-type compartment model. The orange frame marks those agent states that can give positive COVID-19 tests. }
    \label{states}
\end{figure}

\begin{itemize}
    \item Agent’s epidemiological status. Each agent may have one of the 10 stages of the disease (Fig.~\ref{states}) $\vec X = (S, E, I, A, Y, M, H, C, R, D)$.

    \item Agent’s chance to be tested for COVID-19  ($\tilde{p}(X(t))$) that is determined by the agent’s epidemiological status. The agents are tested daily, the number of tests corresponds to the statistical data obtained in a region. At every modeling step, the tests are distributed across the population, the agents whose status is marked with an orange frame in a Fig.~\ref{states} can be given a positive result. The agents whose test is recognized as positive are marked as “confirmed” and included in {\fontfamily{cmss}\selectfont New Diagnoses}. The model assumes that the likelihood for an agent to be tested as a symptomatic carrier is higher and this chance ratio is controlled by parameter $\tilde{p}(X(t))$ at is restored from solving the inverse problem (see the next section). 
    \end{itemize}

\subsection{Parameter identification problem} \label{sect: inv problem}

The ABM developed was also characterized by unknown parameters vector $ \vec q $. To specify the model’s parameters a variational inverse problem formulation was performed to minimize the misfit function: 
\begin{equation}\label{funcJ}
J(\vec q) = \sum_{s}\sum_{t_i=1}^{T}{|X^s_d(t_i) - X^s_m(t_i, \vec q)| \over M_s}.
\end{equation}
Here,  $s$ is the statistics used for comparison, e.g. {\fontfamily{cmss}\selectfont New Diagnoses}, {\fontfamily{cmss}\selectfont New Deaths}, etc.; 
 $X^s_d(t_i)$ and $X^s_m(t_i, \vec q)$ are daily statistical and model data for statistics $s$; $T$ -- the number of modeled days, $M_s = \max\limits_{{t_i}}\{X^s_d(t_i)\}$ is a normalising term.

\subsection{Automatic parameter calibration}

The model assumed that parameter $\beta$ was a piece-wise constant, so the days of change and new values were determined by parameters $\vec \beta_d$ and $\vec\beta_c$. The longer was a considered time interval, the more unknown parameters it included. Since every launch of the model’s calibration algorithm was rather time-consuming, the time interval in question was divided into periods of 1 month. For example, for NY State the first period was 2020.03-02 - 2020.04.01, the second – 2020.04.02 - 2020.05.0; for the UK  – 2020.02.07 - 2020.03.08 and 2020.03.09 - 2020.04.07; for the Novosibirsk region - 2020.03.12 - 2020.04.11 and 2020.04.12 - 2020.05.11 etc.

Each period was sequentially calibrated, so the parameters restored at a previous step were used in the following iteration of the optimization algorithm. Thus, for the initial period considered the unknown parameter vector was
\begin{equation*}
\vec q_1 = (E(0), \beta, \beta_d(1), \beta_c(1), \tilde{p}(X)),
\end{equation*}
where $E(0)$ is the initial number of infected agents, $\beta$ is a contagiousness parameter value, $\beta_d(1)$ is the day parameter $\beta$ changes, $\beta_c(1)$ is the value by which parameter $\beta$ changes on day $\beta_d$, $\tilde{p}(X)$ is a test level parameter in relation to statistical data. For all the following periods (second, third, etc.):
\begin{equation*}
\vec q_i =(\beta_d(i), \beta_c(i)).
\end{equation*}

\begin{figure}[H]
\begin{minipage}[h]{1\linewidth}
\center{\includegraphics[width=\textwidth]{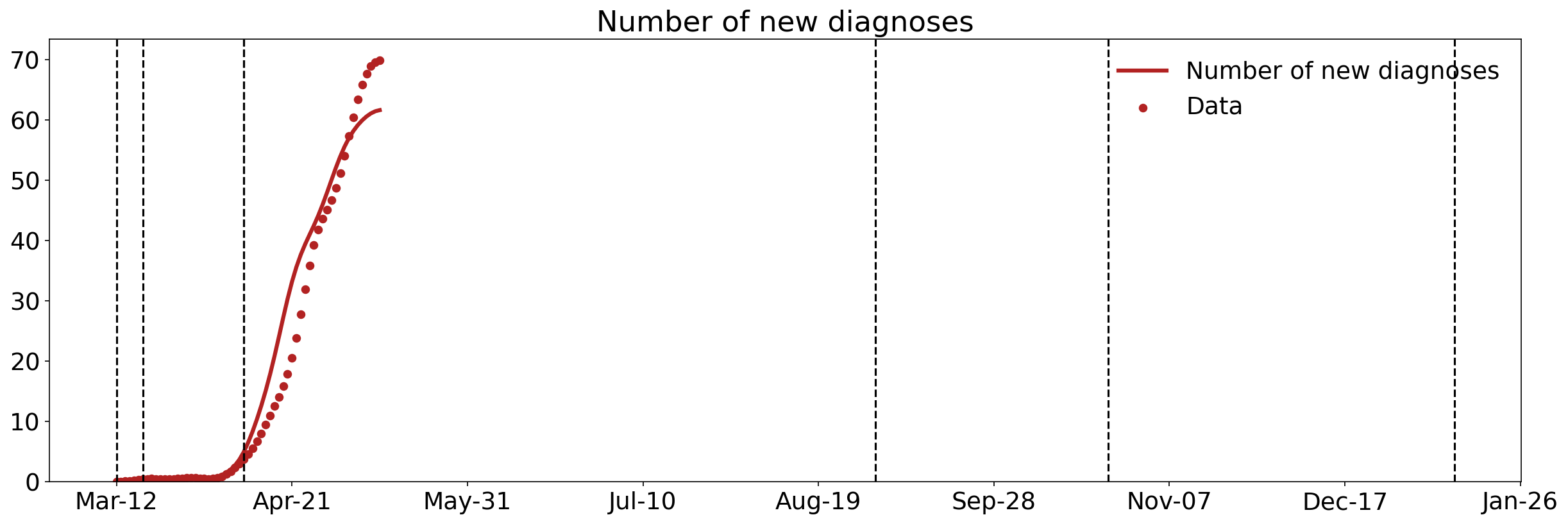}} 
\end{minipage}
\vfill
\begin{minipage}[h]{1\linewidth}
\center{\includegraphics[width=\textwidth]{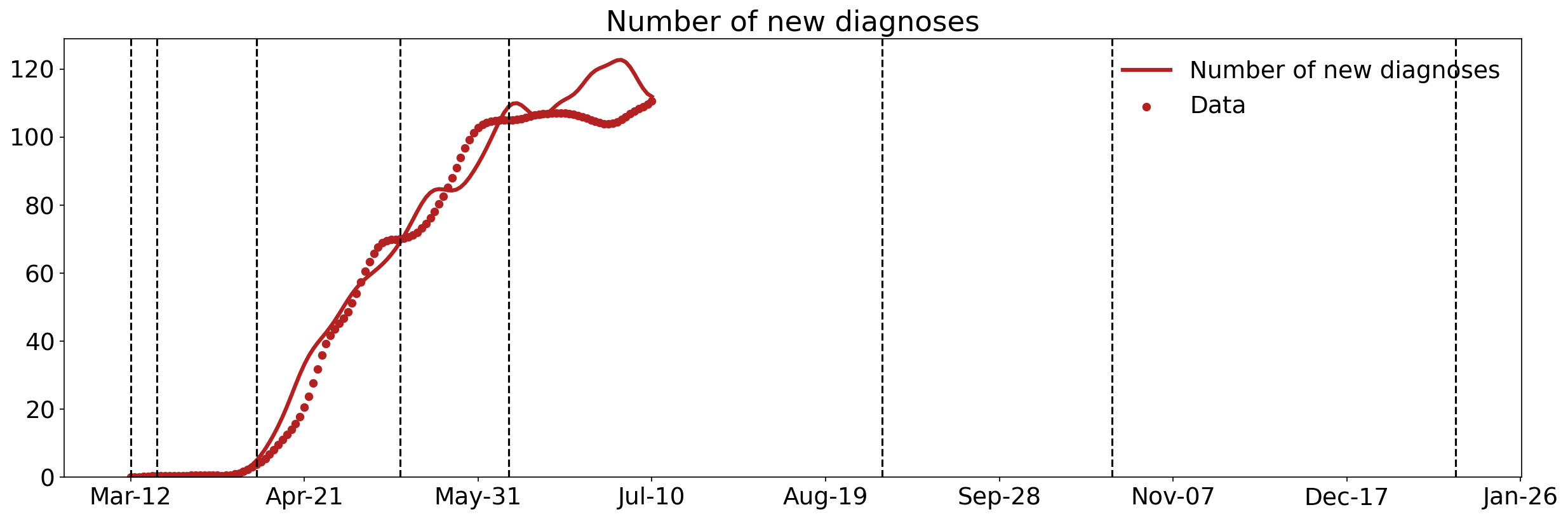}} 
\end{minipage}
\vfill
\begin{minipage}[h]{1\linewidth}
\center{\includegraphics[width=\textwidth]{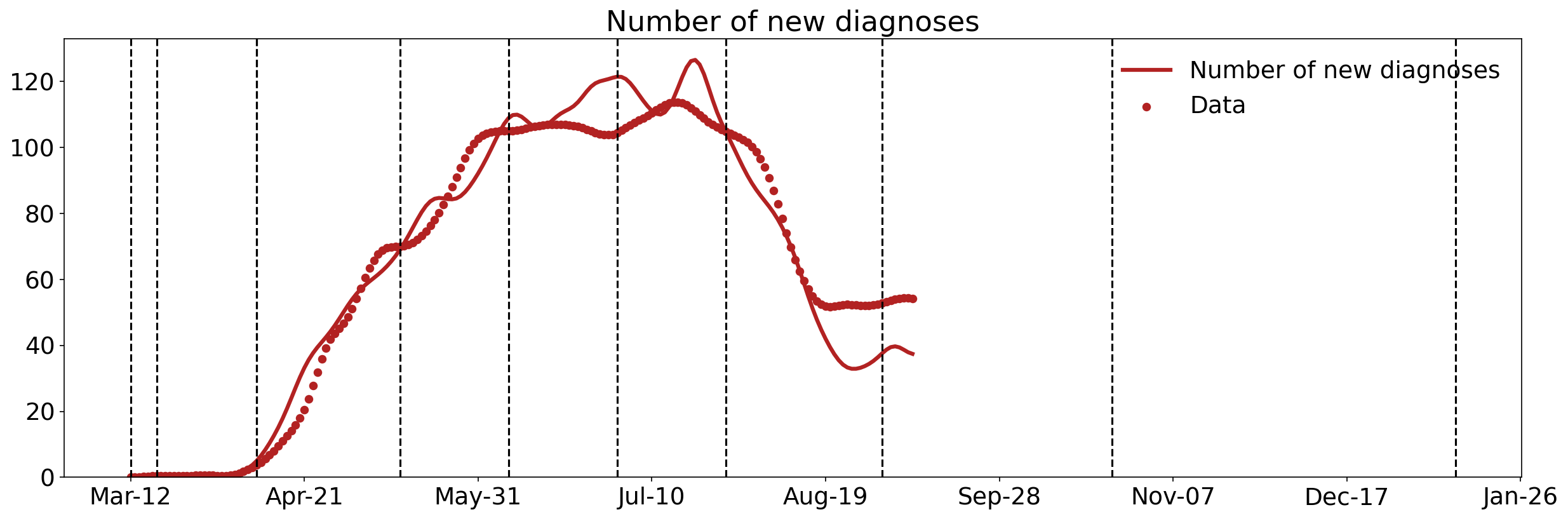}} 
\end{minipage}
\vfill
\begin{minipage}[h]{1\linewidth}
\center{\includegraphics[width=\textwidth]{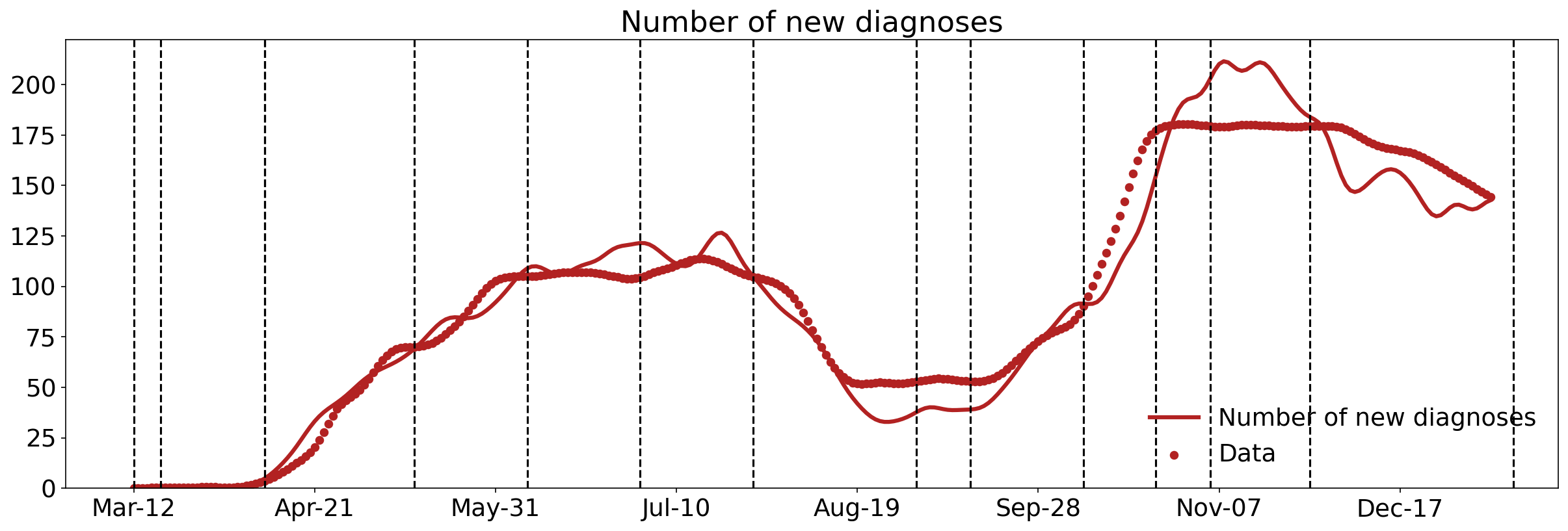}} 
\end{minipage}
\vfill
\caption{Phased restoration of parameters $\vec q$ using the {\fontfamily{cmss}\selectfont New Diagnoses} dataset from the Novosibirsk region. The vertical dashed lines mark the days when the parameter $\beta$ changed that were restored from solving an optimization problem.}
\label{calibration process}
 \end{figure}

Figure~\ref{calibration process} displays a phased restoration of the vector of unknown parameters  $\vec q$ for the Novosibirsk region exemplified by 4 out of 13 interim calibration periods.

%% file: 5methods_and_approaches.tex
\section{Methods and approaches} \label{sect: methods}

Identifying (calibrating) the parameters of an ABM, so that its outputs match observed data, is quite a complex and computation-intensive task due to the large number of the parameters involved. There are different approaches to the problem. In most cases, the parameters are selected manually or one uses averaged experimental results neglecting the specific features of a studied region are used. Analysis of the papers describing epidemiological ABMs has demonstrated that no algorithm could be considered superior for identification of model parameters ~\cite{cal_ABM}. According to the paper “...it appears that calibrating individual-based models in epidemiological studies of HIV, malaria and TB transmission dynamics remains more of an art than a science."

In our model, the vector of unknown parameters $\vec q$ was calibrated using the Optuna hyperparameter optimization software \cite{OPTUNA} to be one of the latest optimizers designed to adjust hyperparameters in machine learning algorithms and neural networks. The optimizer is based on the tree-structured Parzen estimator (TPE) that in many ways is similar to the Bayesian optimizer \cite{Bayes}.
However, unlike the Bayesian optimizer that calculates $p(J(\vec q)|\vec q)$, TPE calculates $p(\vec q|J(\vec q))$ and $p(J(\vec q))$ to determine the parameters domain to minimize functional $J$ by performing Parzen window density estimation, to generate two separate distributions specifying the high - and low-quality regions of the input-space respectively.

The TPE works by ranking the accumulated observations $\mathcal{D}_K = \{q_k
, J(q_k):k = 1, ...,K\}$ based on the objective function values. The algorithm selects where $J$ should be evaluated next by splitting $\mathcal{D}_K$ into two sets, $\mathcal{D}^l_{K_l}$ and $\mathcal{D}^g_{K_g}$ such that $\mathcal{D}^l_{K}$ contains the $\gamma$-percentile of the highest
quality (i.e. lowest function value) points of $\mathcal{D}_K$, whereas $\mathcal{D}^g_{K}$ contains the remaining points. With these definitions, the prior for $J$ is given by $P(J < J^\gamma) = \gamma$. Now, applying non-parametric adaptive parzen windows, two probability distributions $l(\vec q)$ and $g(\vec q)$ are estimated from $\mathcal{D}^l_{K_l}$
and $\mathcal{D}^g_{K_g}$, respectively. Here, $l(\vec q)$ is interpreted as representing
the probability of a region in the input space yielding a high-quality observation while similarly, $g(\vec q)$ represents low-quality regions.
The likelihood of $\vec q_{K+1}$ belonging to distribution $l$ and $g$ may now be expressed by:
\begin{equation}
    p(\vec q_{K+1}|J(\vec q)) = \begin{cases}
                l(\vec q), &  J(\vec q_{K+1}) <J^\gamma,\\
                g(\vec q), & J(\vec q_{K+1}) \geq J^\gamma,
                \end{cases}
\end{equation}
where $J^\gamma$ is the lowest function value found in $\mathcal{D}^g_{K_g}$.

The main optimization criterion in the TPE is an expected improvement (EI).
According to the definition~\cite{TPE}, the EI is a value the function reduces by at a given
moment and can be written as:

\begin{equation}
    EI(\vec q) = \left(\gamma +\frac{g(\vec q)}{l(\vec q)}(1-\gamma)\right)^{-1}
\end{equation}

Finally, in order to select a new point to evaluate the objective function at, samples
are drawn from $l(\vec q)$, where the sample generating the highest expected improvement is
used to evaluate the objective function ($n_{samp}$). A full TPE optimization procedure is described in Algorithm~\ref{TPE}.

\begin{algorithm}
  \caption{Tree-Parzen estimator optimization}\label{TPE}
  \begin{algorithmic}[1]
  \Require{Parameter values for $\gamma$, $n_{samp}$ and $max\_iter$}
  \State \textbf{Inintialize:} accumulate initial observations 
  
     \State $\mathcal{D}_{init}=\{\vec q_k, J(\vec q_k), k=1,\ldots ,n_{init}\}$
  \For{$m$=0 to $max\_iter$} 
    \State Split $\mathcal{D}_{n_{init}+m}$ to generate $\mathcal{D}^g_{m_{g}}, \mathcal{D}^l_{m_l}$
    \State Estimate $l(\vec q)$ from $\mathcal{D}^l_{n_{init}+m_l}$
    \State Estimate $g(\vec q)$ from $\mathcal{D}^g_{n_{init}+m_g}$
    \State Draw $\vec q^s={\vec q^s_k}:k=1,\ldots, n_{samp}$, where $\vec q^s_k \sim l(\vec q)$
    \State $\vec q_{m+1} = \mbox{argmax} EI(\vec q)$
    \State Evaluate $J(\vec q_{m+1})$
    \State Augment set of observations $\mathcal{D}_{n_{init}+m} \leftarrow \mathcal{D}_{n_{init}+m+1} $
  \EndFor
  \end{algorithmic}
\end{algorithm}

%% file: 6modelling_and_forecasting.tex
\section{Modeling and forecasting} \label{sect: results}
In this section, we consider mathematical models and scenarios of COVID-19 spread in NY State, (Section~\ref{sec_NYstate}), the UK (Section \ref{sec_UK}) and the Novosibirsk Region (Section \ref{sec_NSO}) 

\subsection{Initial datasets}

To build and analyze the ABMs in the 3 considered regions, the following data were used:  

\begin{enumerate}
\item Information on population’s age distribution according to the local government statistics;
\item Information on the average family size according to the UN data \cite{UN};
\item Statistical data on the people infected with COVID-19, who recovered and died including the number of tests performed that were collected from:
    \begin{itemize}
        \item The COVID Tracking Project (New York State):
        
        \url{https://covidtracking.com/data}; 
        \item The official UK Government website for data and insights on Coronavirus (United Kingdom):
        
        \url{https://coronavirus.data.gov.uk/};
      
        \item The RBC website and are available for downloading (Novosibirsk region):
        
        \url{http://covid19-modeling.ru/data/novosibirsk-region-data.csv}. 
    \end{itemize}
  
\end{enumerate}

For every region, the modeling results for the  {\fontfamily{cmss}\selectfont New Diagnoses}, {\fontfamily{cmss}\selectfont New Deaths} and  {\fontfamily{cmss}\selectfont Num Critical} (require ICU) datasets were analyzed. In Sections ~\ref{sec_NYstate}-\ref{sec_NSO} one can find the graphs of 30-day forecasts validated with historical data. The forecasts have an 80\% confidence interval to characterize 10\% and 90\% quantiles. In addition, effective reproduction numbers were calculated for every region. These numbers indicate how many persons an infectious agent infects on average during the time it has been infectious and calculated as: 
\begin{equation}\label{R_eff_formula}
R(t)=\frac{I_N(t)\cdot f}{I_C(t)},
\end{equation}
where $I_N(t)$ is the number of new infections on day $t$, $I_C(t)$ is the number of actively infectious people on day $t$ and $f$ is the average duration of infectiousness. If $R(t)<1$, the pandemic is considered to stop spreading and keeps spreading otherwise.


\subsection{COVID-19 spread simulation in NY State}\label{sec_NYstate}
The mismatch (also called error or loss) function~(\ref{funcJ}) for New York State was as follows:
\begin{equation*}
    J(\vec q) = \sum_{t_i=1}^{T}\left(\frac{|Y_d(t_i) - Y_m(t_i, \vec q)|}{M_{diag}} + \frac{|D_d(t_i) - D_m(t_i, \vec q)|}{M_{death}}\right).
\end{equation*}

Here, $Y_d(t_i), Y_m(t_i, \vec q)$ are smoothed {\fontfamily{cmss}\selectfont New Diagnoses} with confirmed COVID-19, $D_d(t_i), D_m(t_i, \vec q)$ are smoothed {\fontfamily{cmss}\selectfont New Deaths}.

\begin{figure}[H]
    \centering
    \includegraphics[width=\textwidth]{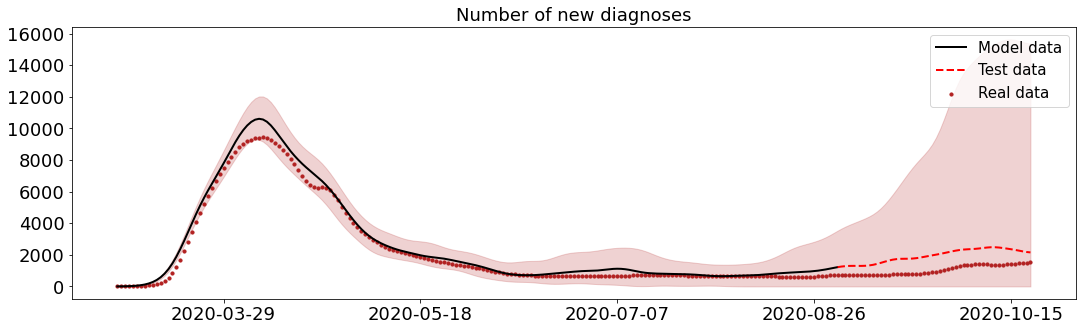}
    \caption{
     Forecast for {\fontfamily{cmss}\selectfont New Diagnoses} validated with historical data for NY State. The black line displays the modeled period, and the red dots – the {\fontfamily{cmss}\selectfont New Diagnoses} dataset.}
    \label{prog_hist_NY}
\end{figure}

After identification of the parameters, the model was validated with historical data. The {\fontfamily{cmss}\selectfont New Diagnoses} dataset for NY State from 2020.03.02 to 2020.08.31 was taken as model data. The produced forecast covered 50 days and was compared against statistical data for the given test period (2020.09.01 –- 2020.10.20). The modeling results can be seen in Fig.~\ref{prog_hist_NY}. 

Despite the size of the confidence interval that started to increase from 2020.09.01, the increment rate of New Diagnoses was close to the test data.

\begin{figure}[H]
\begin{minipage}[h]{1\linewidth}
\center{\includegraphics[width=\textwidth]{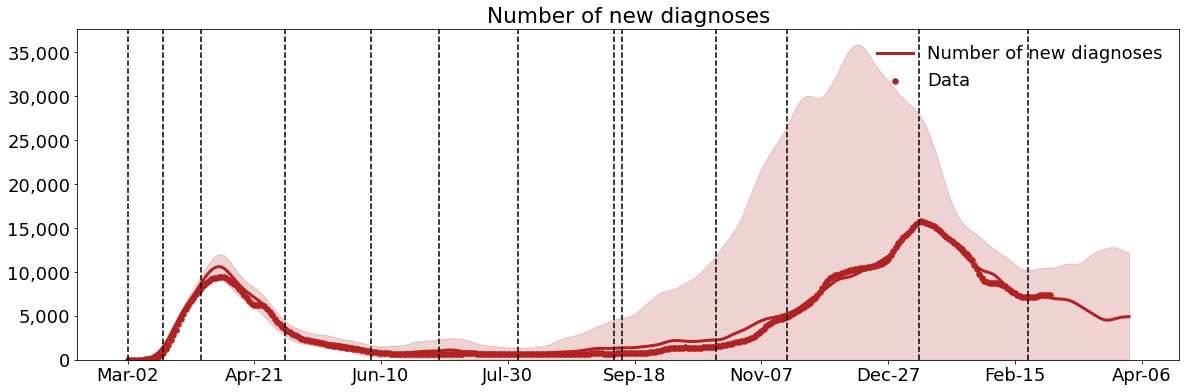}} 
\end{minipage}
\vfill
\begin{minipage}[h]{1\linewidth}
\center{\includegraphics[width=\textwidth]{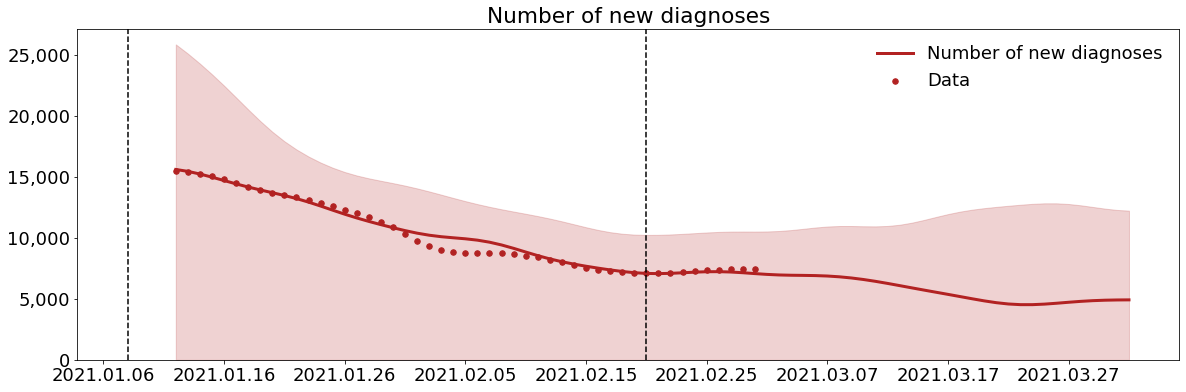}} 
\end{minipage}
\vfill
\begin{minipage}[h]{1\linewidth}
\center{\includegraphics[width=\textwidth]{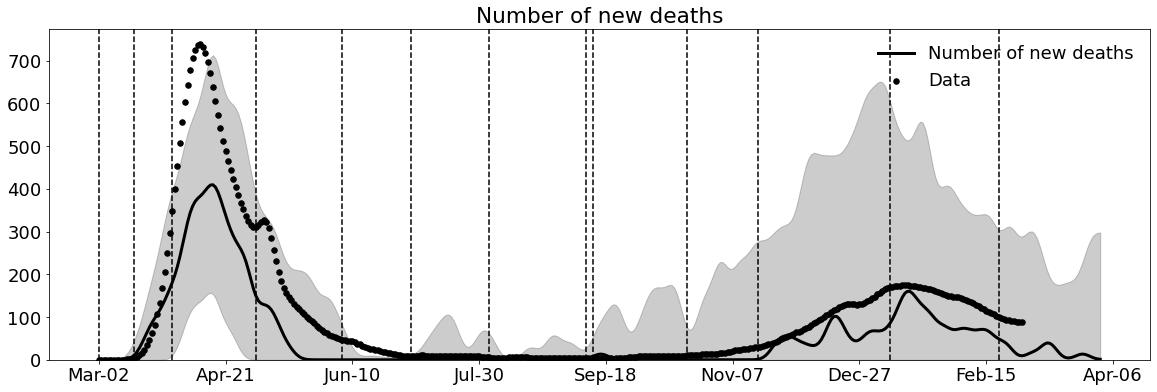}} 
\end{minipage}
\vfill
\begin{minipage}[h]{1\linewidth}
\center{\includegraphics[width=\textwidth]{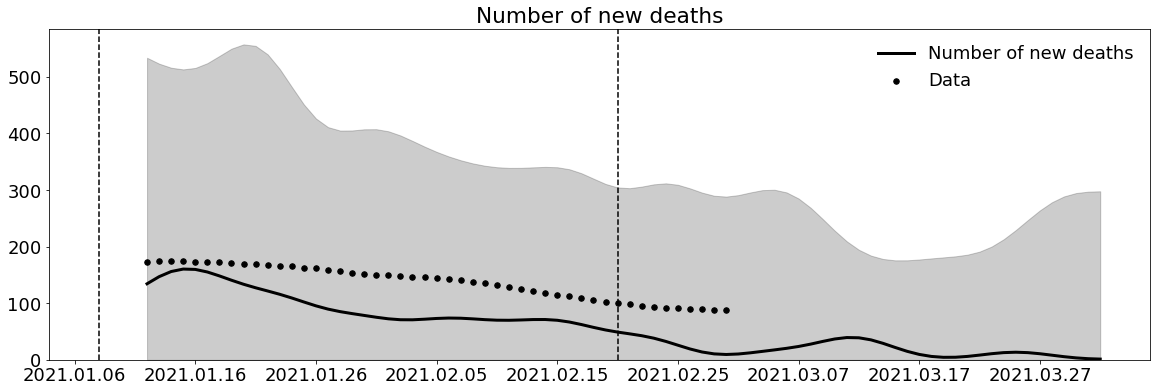}} 
\end{minipage}
\vfill
\caption{Model calibration results for 10 simulations + a 30-day forecast for {\fontfamily{cmss}\selectfont New Diagnoses} (two top graphs) and {\fontfamily{cmss}\selectfont New Deaths} (two bottom graphs) in NY State. The shaded areas are 10\% and 90\% quantiles, the solid line -- the median value of modeling result, and dots -- real data. The dashed vertical lines are COVID-19 containment measures.}
\label{model_plot_NY}
 \end{figure}

Figure~\ref{model_plot_NY} presents the results of a 30-day forecast with restored vector of unknown parameters for NY State, where the dots mark the real data accumulated from 2020.03.02 to 2021.03.01. The forecast for {\fontfamily{cmss}\selectfont New Diagnoses} assumed the rate of daily tests remained unchanged.

The effective reproduction number  $R(t)$ calculated by~(\ref{R_eff_formula}) showed the pandemic began its downturn in NY State in January 2021. 

\begin{figure}[h!]
    \centering
    \includegraphics[width=\textwidth]{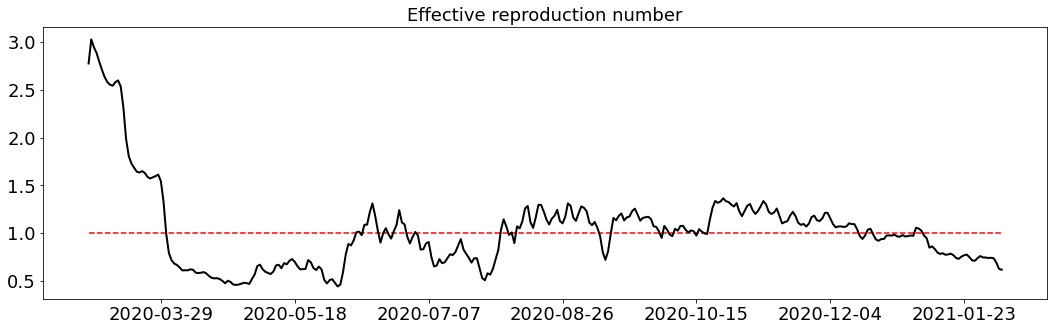}
    \caption{Effective reproduction number $R(t)$ for the New York State. The red dashed line marks $R(t)=1$.}
    \label{r_eff_NY}
\end{figure}


\subsection{COVID-19 spread simulation in the UK}\label{sec_UK}
The mismatch function~(\ref{funcJ}) for United Kingdom was as follows:
\begin{equation*}
    J(\vec q) = \sum_{t_i=1}^{T}\frac{|Y_d(t_i) - Y_m(t_i, \vec q)|}{M_{diag}}.
\end{equation*}
Here $Y_d(t_i), Y_m(t_i, \vec q)$ are smoothed {\fontfamily{cmss}\selectfont New Diagnoses} with confirmed COVID-19. 

\begin{figure}[H]
    \centering
    \includegraphics[width=\textwidth]{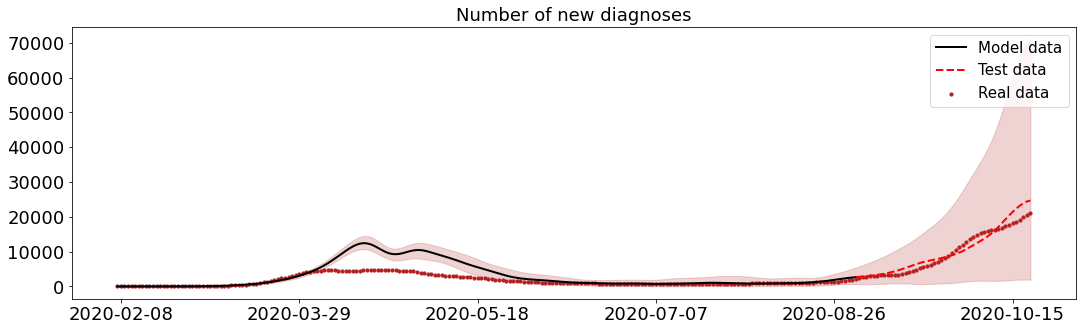}
    \caption{Forecast for {\fontfamily{cmss}\selectfont New Diagnoses} validated with historical data for UK. The black line displays the modeled period, and the red dots – the {\fontfamily{cmss}\selectfont New Diagnoses} dataset. The red dashed line marks the test period.}
    \label{prog_hist_uk}
\end{figure}

Firtly, the parameters that minimizing mismatch function were identified.
The {\fontfamily{cmss}\selectfont New Diagnoses} dataset for the UK from 2020.02.07 to 2020.08.31 was taken as model data. The produced forecast covered 50 days and was compared against observed data for the given test period (2020.09.01 – 2020.10.20). The modeling results can be seen in Fig.~\ref{prog_hist_uk}. It is noteworthy that the prediction error after comparing with real data was only 1.2\%. 

\begin{figure}[H]
\begin{minipage}[h]{1\linewidth}
\center{\includegraphics[width=\textwidth]{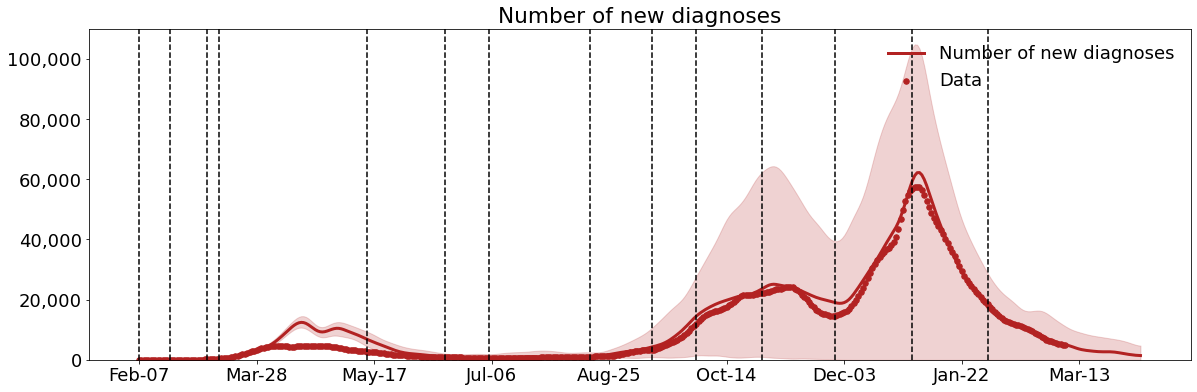}} 
\end{minipage}
\vfill
\begin{minipage}[h]{1\linewidth}
\center{\includegraphics[width=\textwidth]{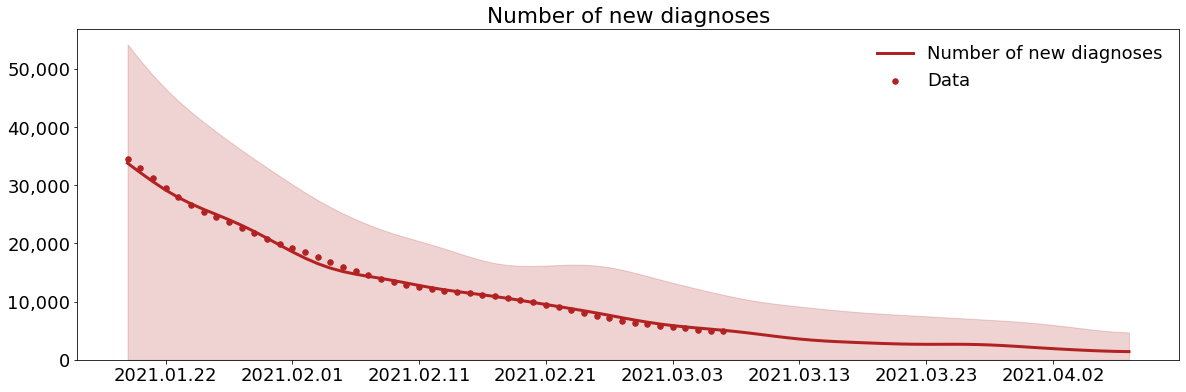}} 
\end{minipage}
\vfill
\begin{minipage}[h]{1\linewidth}
\center{\includegraphics[width=\textwidth]{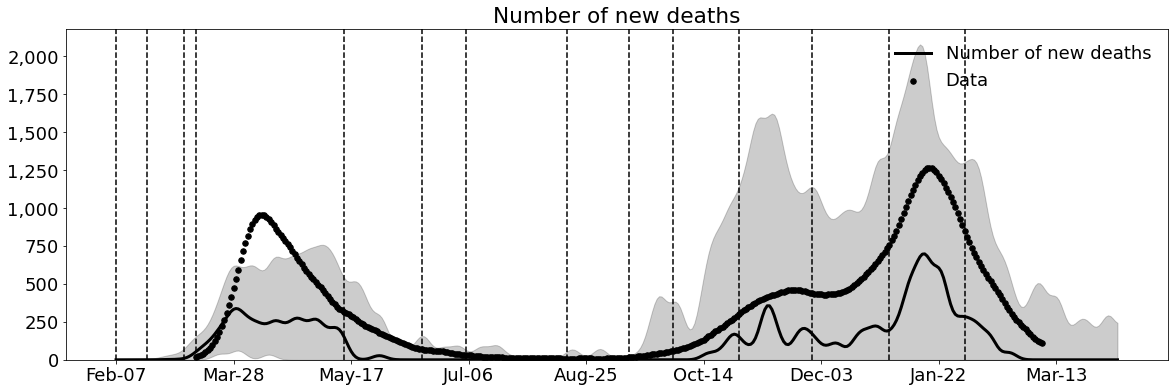}} 
\end{minipage}
\vfill
\begin{minipage}[h]{1\linewidth}
\center{\includegraphics[width=\textwidth]{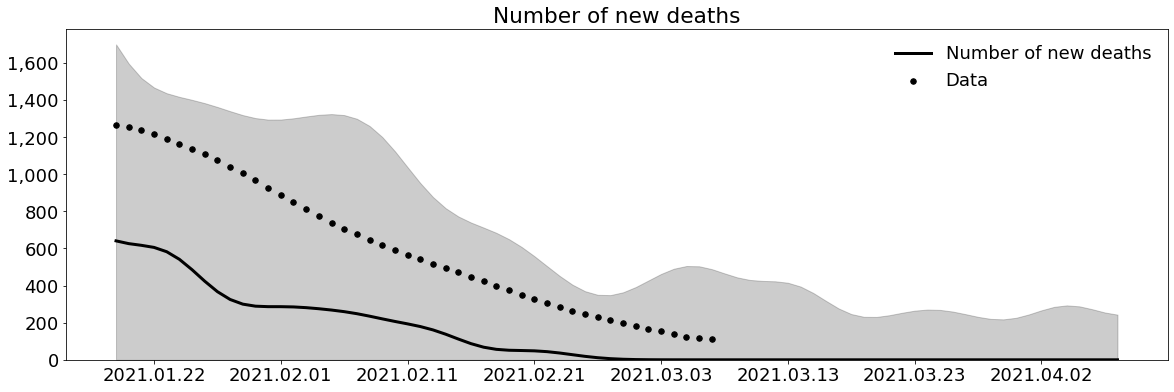}} 
\end{minipage}
\vfill
\caption{Model calibration results for 10 simulations + a 30-day forecast for {\fontfamily{cmss}\selectfont New Diagnoses} and {\fontfamily{cmss}\selectfont New Deaths} in the UK. The shaded areas are 10\% and 90\% quantiles, the solid line -- the median value of modeling result, and dots -- real data. The dashed vertical lines are COVID-19 containment measures.}
\label{model_plot_uk}
 \end{figure}

Figure.~\ref{model_plot_uk} presents the results of a 30-day forecast for the UK with restored vector of unknown parameters $\vec q$ as is described in Section 4.  The forecast for {\fontfamily{cmss}\selectfont New Diagnoses} assumed the rate of daily tests remained unchanged.

As for {\fontfamily{cmss}\selectfont New Deaths}, the modeling results differed from the real data by 50\%, which was related to the target functional $J(\vec q)$ that contained information only about {\fontfamily{cmss}\selectfont New Diagnoses}, while the error for this indicator was comparable to that of {\fontfamily{cmss}\selectfont New Deaths}. To reduce the error, weighed information about the {\fontfamily{cmss}\selectfont New Deaths} indicator is to be added in the functional. The problem of determining the weight coefficients needs further consideration.

\begin{figure}[h!]
    \centering
    \includegraphics[width=\textwidth]{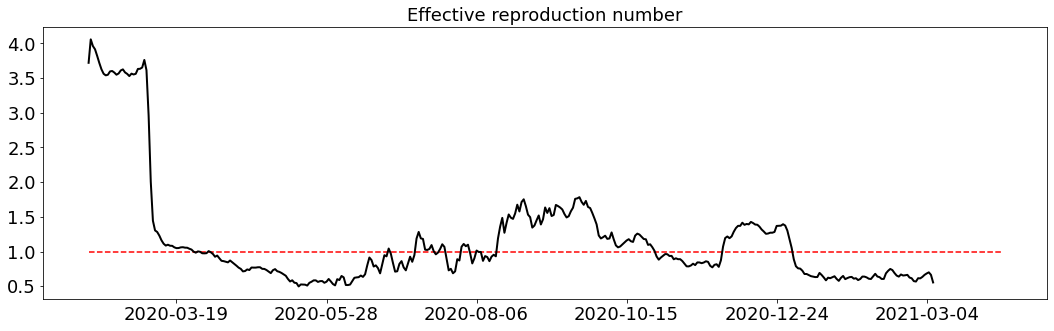}
    \caption{Effective reproduction number $R(t)$ for the UK. The red dashed line marks $R(t)=1$.}
    \label{r_eff_UK}
\end{figure}

The effective reproduction number $R(t)$ calculated by (\ref{R_eff_formula}) for a time period from 2020.02.07 to 2021.03.01 in the UK can be seen in Fig.~\ref{r_eff_UK}. The graph demonstrates the UK had a second infection wave in fall 2020 and in January 2021, after which the pandemic began its downturn.

\subsection{COVID-19 spread simulation in the Novosibirsk region}\label{sec_NSO}
The mismatch function~(\ref{funcJ}) for Novosibirsk region was as follows:
\begin{equation*}
    J(\vec q) = \sum_{t_i=1}^{T}\frac{|Y_d(t_i) - Y_m(t_i, \vec q)|}{M_{diag}}.
\end{equation*}
Here $Y_d(t_i), Y_m(t_i, \vec q)$ are smoothed {\fontfamily{cmss}\selectfont New Diagnoses} with confirmed COVID-19. 

\begin{figure}[H]
    \centering
    \includegraphics[width=\textwidth]{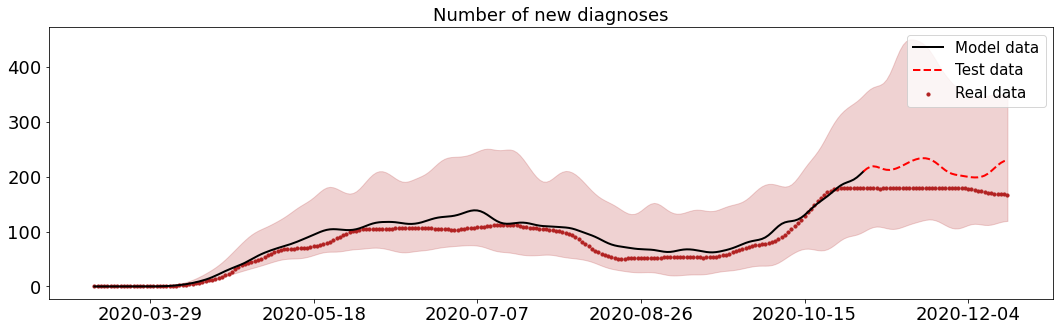}
    \caption{
    Forecast for {\fontfamily{cmss}\selectfont New Diagnoses} validated with historical data for the Novosibirsk region. The black line displays the modeled period, and the red dashed line – the test period. The red dashed line marks the test period. The shaded areas are 10\% and 90\% quantiles for 10 simulations}
    \label{prog_hist_Novosib}
\end{figure}

The {\fontfamily{cmss}\selectfont New Diagnoses} dataset for validation model in the Novosibirsk region was from 2020.03.12 to 2020.11.01. The produced forecast covered 45 days and was compared against statistical data for the given test period (2020.11.02 – 2020.12.15). The modeling results shown in Fig.~\ref{prog_hist_Novosib} demonstrated that preserving the containment measures introduced in the region could have reduced the spread by more than 25\%. This difference between modeling and real data was due to the data space’s stationarity (see Section 2.2.2) amid unaccounted administrative measures. 

\begin{figure}[H]
\begin{minipage}[h]{1\linewidth}
\center{\includegraphics[width=\textwidth]{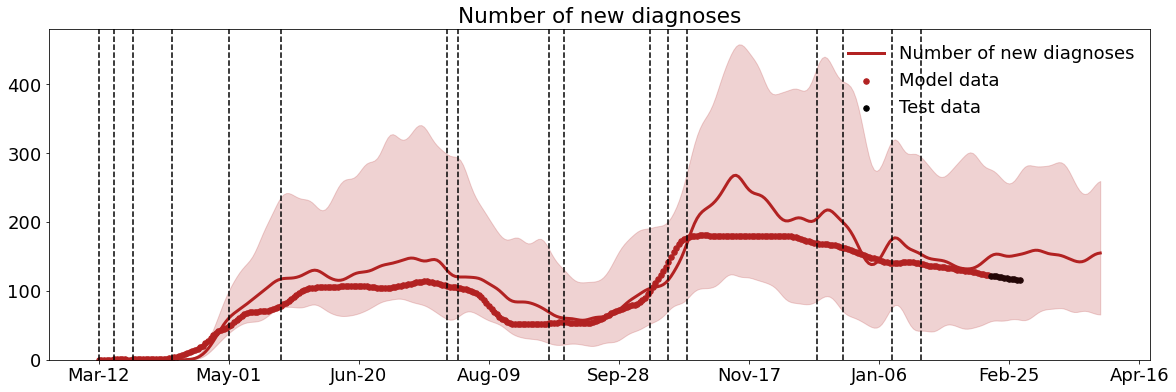}} 
\end{minipage}
\vfill
\begin{minipage}[h]{1\linewidth}
\center{\includegraphics[width=\textwidth]{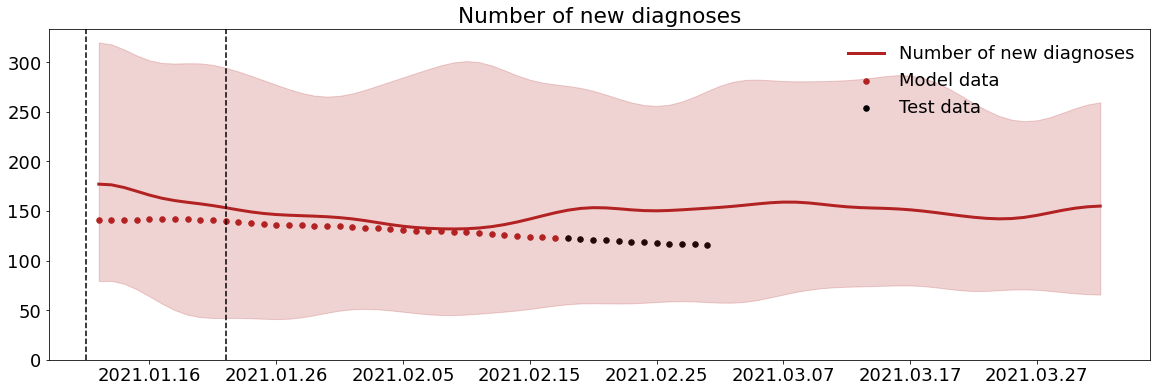}} 
\end{minipage}
\vfill
\begin{minipage}[h]{1\linewidth}
\center{\includegraphics[width=\textwidth]{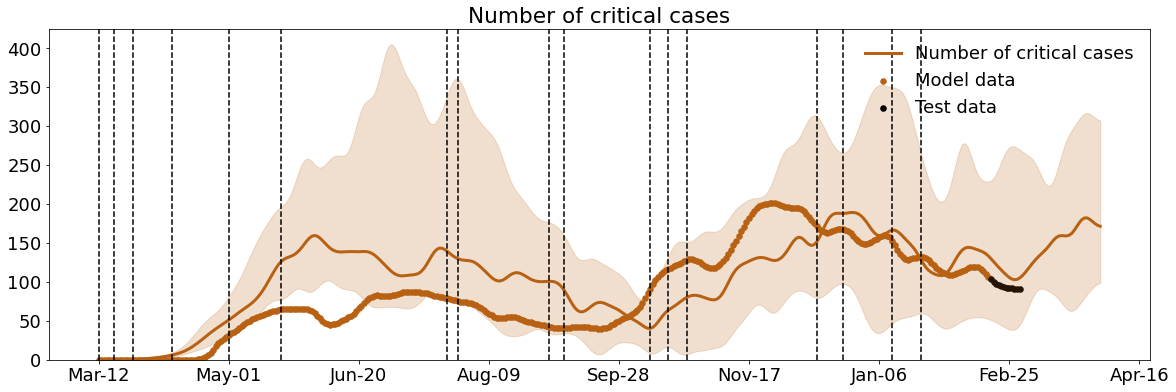}} 
\end{minipage}
\vfill
\begin{minipage}[h]{1\linewidth}
\center{\includegraphics[width=\textwidth]{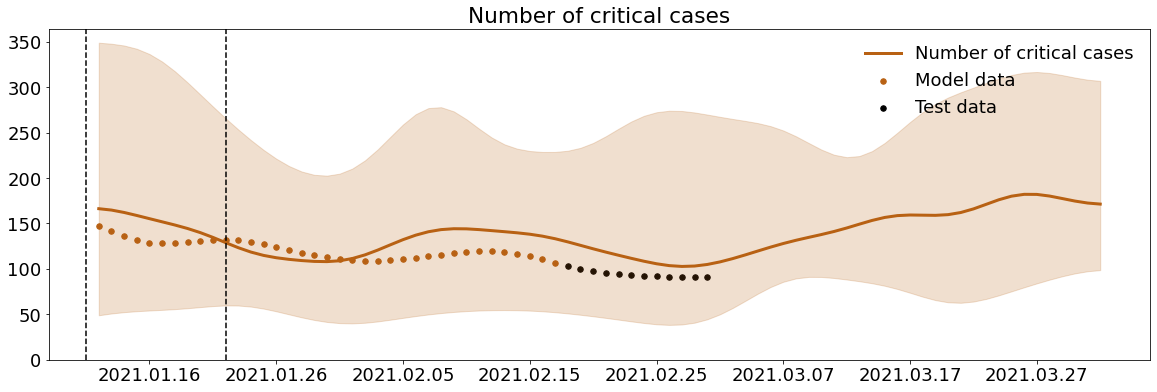}} 
\end{minipage}
\vfill
\caption{Model calibration results for 10 simulations + a 30-day forecast for {\fontfamily{cmss}\selectfont New Diagnoses} and {\fontfamily{cmss}\selectfont Num Critical} in the Novosibirsk region. The shaded areas are 10\% and 90\% quantiles, the solid line -- the median value of modeling result, and dots -- real data. The dashed vertical lines are COVID-19 containment measures.}
\label{model_plot_Novosibirsk}
 \end{figure}

Figure~\ref{model_plot_Novosibirsk} presents the results of a 30-day forecast for the Novosibirsk region with a restored vector of unknown parameters $\vec q$ as is described in Section 4. The colored dots mark the real data (to 2021.02.15) that were used to solve the inverse problem; the black ones (2021.02.16 – 2021.03.01) – to test the forecast. The forecast for {\fontfamily{acmss}\selectfont New Diagnoses} assumed the rate of daily tests remained unchanged. It is noteworthy that despite the forecast error for the {\fontfamily{cmss}\selectfont New Diagnoses} indicator participating in the minimization of target functional $J(\vec q)$ was 10\%, the forecast error for {\fontfamily{acmss}\selectfont Num Critical} did not exceed 5\%. 

\begin{figure}[h!]
    \centering
    \includegraphics[width=\textwidth]{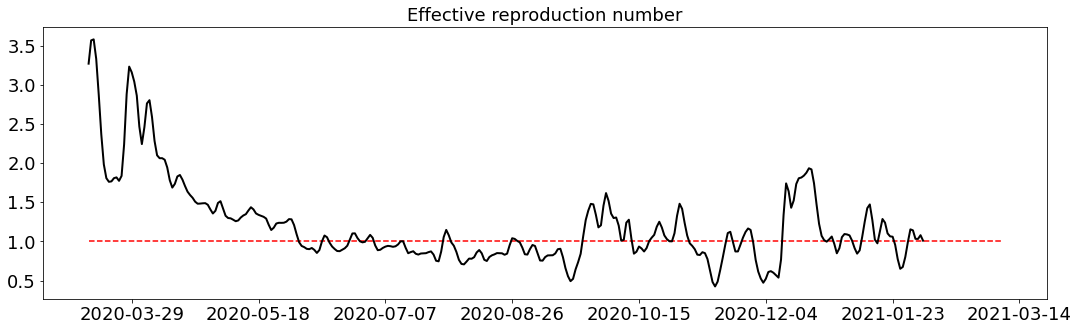}
    \caption{Effective reproduction number $R(t)$ for the Novosibirsk region. The red dashed line marks $R(t)=1$.}
    \label{r_eff_Novosib}
\end{figure}

 The effective reproduction number $R(t)$ calculated by (\ref{R_eff_formula}) for a time period from 2020.03.12 to 2021.02.01 in the Novosibirsk region can be seen in Fig. \ref{r_eff_Novosib}. The graph demonstrates the pandemic had been under control since the end of 2020.

%% file: 7conclusion.tex
\section{Conclusions and discussions}

Mathematical models are effective tools to deal with the time evolution and patterns of disease outbreaks. They provide us with useful predictions in the context of the impact of intervention in decreasing the number of infections and deaths.

This paper provides a detailed analysis of the statistical data, including the number of tested, positive, mortality, hospitalization and critical cases, on COVID-19 spread in NY State, the UK and the Novosibirsk region and presents a solution to the problem of identification of unknown epidemiological parameters (transmissibility; the initial number of infected individuals; probability of being tested, etc.) in an ABM. The problem has been considered as the minimization of a target functional in relation to daily numbers of tested, positive and mortality cases in the studied regions and become an important modification of the Covasim package [21]. The minimization problem has been solved using the Optuna hyperparameter optimization software and the Parzen estimation method. Normal gradient descent methods do not work with Covasim or other agent-based models, due to the stochastic variability between model runs that makes the landscape very “bumpy” (i.e., many transient local minima). One way of getting around this is to perform many different runs and take the average. However, averaging over many runs is computationally expensive, since running $N$ simulations of ABM will only reduce the noise by $\sqrt{N}$ ~\cite{covasim_doc}.

The results of data analysis in every studied region showed the weekly seasonality of {\fontfamily{cmss}\selectfont New Tests}, which helped us forecast the future values of this time series. It worth noting that some countries have their own rules of statistical analysis that have to be accounted for when carrying out modeling, e.g. in the USA, the {\fontfamily{cmss}\selectfont New Diagnoses} indicator contains a certain percentage of probable cases, while {\fontfamily{cmss}\selectfont New Deaths} in the UK accounts for all the deaths that have occurred within 28 days since a positive COVID-19 test, even if such death has not been provoked by the virus. 

Due to the high sensitivity of the transmissibility parameter $\beta$, the accuracy of its identification becomes crucial for uncovering the pattern of COVID-19 spread in an investigated region. However, not all containment measures affect the pattern. For that reason, when developing an ABM calibration algorithm based on epidemiological data, we paid special attention to transmissibility and the times this parameter changed while modeling. These characteristics were determined from solving the minimization problem as a piecewise - constant function, while solving the inverse problem restored the parameter together with its times of change. 

The devised ABM has been validated with historical data from 2020. The modeling results for the three regions in question have demonstrated that preserving the introduced containment measures would have sustained {\fontfamily{cmss}\selectfont New Diagnoses} in NY State and the Novosibirsk region during March 2021 and would have reduced them in the UK.  

Due to the specific features of the observed data in the Novosibirsk region, where two time series are stationary with probability of 1, the model’s forecast accuracy has been lower for new cases of infection and more accurate for cases of hospitalization, e.g., the confidence interval values for the new case of infection registered daily was [70; 270] individuals and for hospitalizations - [90; 300] individuals, both including as model solutions as real data.

The proposed agent-based model has the following limitations. We do not fix population number changing during model year,
consider waning of immunity to coronavirus, or the possibility of re-infection. We also use more simplistic contact structures than in real life. 

Our future plans are investigation of model identifiability to real data and sensitivity analysis. It will also be necessary to investigate the influence of vaccination on COVID-19 propagation.

%% file: 9competing.tex
\section{Competing interests statement}
The authors declare that they have no known competing financial interests or personal relationships that could have appeared to influence the work reported in this paper.

%% file: acknowledges.tex
\section{Acknowledgements}
The data analysis part (section~2) is supported by the Russian Foundation for Basic Research and Royal Society (project no.~21-51-10003). The agent-based mathematical model construction and analysis of numerical results (sections~3, 4, 5) was supported by the Russian Science Foundation (project no.~18-71-10044) and the Royal Society $IEC\backslash R2\backslash 202020$ -- International Exchanges 2020 Cost Share between UK and Russia. Authors would like to thank Professor Sergey Kabanikhin and Doctor Alexey Romanyukha for the problem analysis and fruitful discussions. Also we would like to thank Yan Reznichenko for language help.

%% file: references.tex
%
%
%
%